\newif\ifpagetitre            \pagetitretrue
\newtoks\hautpagetitre        \hautpagetitre={\hfil}
\newtoks\baspagetitre         \baspagetitre={\hfil}
\newtoks\auteurcourant        \auteurcourant={\hfil}
\newtoks\titrecourant         \titrecourant={\hfil}

\newtoks\hautpagegauche       \newtoks\hautpagedroite
\hautpagegauche={\hfil\the\auteurcourant\hfil}
\hautpagedroite={\hfil\the\titrecourant\hfil}

\newtoks\baspagegauche \baspagegauche={\hfil\tenrm\folio\hfil}
\newtoks\baspagedroite \baspagedroite={\hfil\tenrm\folio\hfil}

\headline={\ifpagetitre\the\hautpagetitre
\else\ifodd\pageno\the\hautpagedroite
\else\the\hautpagegauche\fi\fi}

\footline={\ifpagetitre\the\baspagetitre
\global\pagetitrefalse
\else\ifodd\pageno\the\baspagedroite
\else\the\baspagegauche\fi\fi}

\vsize=9.0in\voffset=1cm
\looseness=2


\message{fonts,}

\font\tenrm=cmr10
\font\ninerm=cmr9
\font\eightrm=cmr8
\font\teni=cmmi10
\font\ninei=cmmi9
\font\eighti=cmmi8
\font\ninesy=cmsy9
\font\tensy=cmsy10
\font\eightsy=cmsy8
\font\tenbf=cmbx10
\font\ninebf=cmbx9
\font\tentt=cmtt10
\font\ninett=cmtt9

\font\ninesl=cmsl9
\font\eightsl=cmsl8

\font\nineit=cmti9
\font\eightit=cmti8

\skewchar\ninei='177 \skewchar\eighti='177
\skewchar\ninesy='60 \skewchar\eightsy='60

\def\eightpoint{\def\rm{\fam0\eightrm} 
\normalbaselineskip=9pt
\normallineskiplimit=-1pt
\normallineskip=0pt

\textfont0=\eightrm \scriptfont0=\sevenrm \scriptscriptfont0=\fiverm
\textfont1=\ninei \scriptfont1=\seveni \scriptscriptfont1=\fivei
\textfont2=\ninesy \scriptfont2=\sevensy \scriptscriptfont2=\fivesy
\textfont3=\tenex \scriptfont3=\tenex \scriptscriptfont3=\tenex
\textfont\itfam=\eightit  \def\it{\fam\itfam\eightit} 
\textfont\slfam=\eightsl \def\sl{\fam\slfam\eightsl} 

\setbox\strutbox=\hbox{\vrule height6pt depth2pt width0pt}%
\normalbaselines \rm}

\def\ninepoint{\def\rm{\fam0\ninerm} 
\textfont0=\ninerm \scriptfont0=\sevenrm \scriptscriptfont0=\fiverm
\textfont1=\ninei \scriptfont1=\seveni \scriptscriptfont1=\fivei
\textfont2=\ninesy \scriptfont2=\sevensy \scriptscriptfont2=\fivesy
\textfont3=\tenex \scriptfont3=\tenex \scriptscriptfont3=\tenex
\textfont\itfam=\nineit  \def\it{\fam\itfam\nineit} 
\textfont\slfam=\ninesl \def\sl{\fam\slfam\ninesl} 
\textfont\bffam=\ninebf \scriptfont\bffam=\sevenbf
\scriptscriptfont\bffam=\fivebf \def\bf{\fam\bffam\ninebf} 
\textfont\ttfam=\ninett \def\tt{\fam\ttfam\ninett} 

\normalbaselineskip=11pt
\setbox\strutbox=\hbox{\vrule height8pt depth3pt width0pt}%
\let \smc=\sevenrm \let\big=\ninebig \normalbaselines
\parindent=1em
\rm}

\def\tenpoint{\def\rm{\fam0\tenrm} 
\textfont0=\tenrm \scriptfont0=\ninerm \scriptscriptfont0=\fiverm
\textfont1=\teni \scriptfont1=\seveni \scriptscriptfont1=\fivei
\textfont2=\tensy \scriptfont2=\sevensy \scriptscriptfont2=\fivesy
\textfont3=\tenex \scriptfont3=\tenex \scriptscriptfont3=\tenex
\textfont\itfam=\nineit  \def\it{\fam\itfam\nineit} 
\textfont\slfam=\ninesl \def\sl{\fam\slfam\ninesl} 
\textfont\bffam=\ninebf \scriptfont\bffam=\sevenbf
\scriptscriptfont\bffam=\fivebf \def\bf{\fam\bffam\tenbf} 
\textfont\ttfam=\tentt \def\tt{\fam\ttfam\tentt} 

\normalbaselineskip=11pt
\setbox\strutbox=\hbox{\vrule height8pt depth3pt width0pt}%
\let \smc=\sevenrm \let\big=\ninebig \normalbaselines
\parindent=1em
\rm}

\message{fin format jgr}

\hautpagegauche={\hfill\ninerm\the\auteurcourant}
\hautpagedroite={\ninerm\the\titrecourant\hfill}
\auteurcourant={R.G.Novikov}
\titrecourant={On non-overdetermined inverse scattering at zero energy
in three dimensions}

\magnification=1200
\font\Bbb=msbm10
\def\C{\hbox{\Bbb C}}
\def\R{\hbox{\Bbb R}}
\def\S{\hbox{\Bbb S}}
\def\N{\hbox{\Bbb N}}
\def\b{\backslash}
\def\v{\varphi}
\def\pa{\partial}
\def\ep{\varepsilon}
\def\pr{\parallel}

\vskip 2 mm
\centerline{\bf On  non-overdetermined inverse scattering}
\centerline{\bf at zero energy in three dimensions}
\vskip 2 mm
\centerline{\bf R. G. Novikov}
\vskip 2 mm
\centerline{31 May 2006}

\noindent
{\ninerm CNRS, Laboratoire de Math\'ematiques Jean Leray (UMR 6629),
Universit\'e de Nantes, BP 92208,}

\noindent
{\ninerm F-44322, Nantes cedex 03, France}

\noindent
{\ninerm e-mail: novikov@math.univ-nantes.fr}
\vskip 4 mm
{\bf Abstract.}
We develop the $\bar\pa$- approach to inverse scattering at zero energy
in dimensions $d\ge 3$ of [Beals, Coifman 1985], [Henkin, Novikov 1987] and
[Novikov 2002]. As a result we give, in particular, uniqueness theorem,
precise reconstruction procedure, stability estimate and approximate
reconstruction for the problem of
finding a sufficiently small potential $v$ in the Schr\"odinger equation
from a fixed non-overdetermined  ("backscattering" type)
restriction $h\big|_{\Gamma}$ of the Faddeev
generalized scattering amplitude $h$ in the complex domain at zero energy
in dimension $d=3$. For sufficiently small potentials $v$ we formulate also
a characterization theorem for the aforementioned restriction
$h\big|_{\Gamma}$ and a new characterization theorem for the full Faddeev
function
$h$ in the complex domain at zero energy in dimension $d=3$. We show
that the results of the present work have direct applications to the
electrical impedance tomography via a reduction given first in [Novikov, 1987,
1988].

\vskip 4 mm
{\bf 1.Introduction}

Consider the Schr\"odinger equation at zero energy
$$-\Delta\psi+v(x)\psi=0,\ \ x\in\R^d,\ \ d\ge 2,\eqno(1.1)$$
where
$$\eqalign{
&v\ \ {\rm is\ a\ sufficiently\ regular\ function\ on}\ \ \R^d \cr
&{\rm with\ sufficient\ decay\ at\ infinity} \cr}\eqno(1.2)$$
(precise assumptions on $v$ are specified below in this introduction and
in Sections 2 and 3). For equation (1.1), under assumptions (1.2), we consider the
Faddeev generalized scattering amplitude $h(k,l)$, where $(k,l)\in\Theta$,
$$\Theta=\{k\in\C^d,\ l\in\C^d:\ k^2=l^2=0,\ \ Im\,k=Im\,l\}.\eqno(1.3)$$
For definitions of $h$ see, for example, [HN] (Section 2.2) and [No1]
(Section 2). Given $v$, to determine $h$ on $\Theta$ one can use, in
particular, the formula
$$h(k,l)=H(k,k-l),\ \ (k,l)\in\Theta,\eqno(1.4)$$
and the linear integral equation
$$H(k,p)=\hat v(p)-\int\limits_{\R^d}
{\hat v(p+\xi)H(k,-\xi)d\xi\over {\xi^2+2k\xi}},\ \ k\in\Sigma,\ \
p\in\R^d,\eqno(1.5)$$
where
$$\eqalignno{
&\hat v(p)=(2\pi)^{-d}\int\limits_{\R^d}e^{ipx}u(x)dx,\ \ p\in\R^d,&(1.6)\cr
&\Sigma=\{k\in\C^d:\ k^2=0\}.&(1.7)\cr}$$

In the present work we consider, mainly, the three dimensional case $d=3$.
In addition, in the main considerations of the present work for $d=3$ our
basic assumption on $v$ consists in the following condition on its Fourier
transform
$$\hat v\in L_{\mu}^{\infty}(\R^3)\ \ {\rm for\ some\ real}\ \ \mu\ge 2,
\eqno(1.8)$$
where
$$\eqalign{
&L_{\mu}^{\infty}(\R^d)=\{u\in L^{\infty}(\R^d):\ \|u\|_{\mu}<+\infty\},\cr
&\|u\|_{\mu}=ess\,\sup\limits_{p\in\R^d}(1+|p|)^{\mu}|u(p)|,\ \ \mu>0.\cr}
\eqno(1.9)$$
If $v$ satisfies (1.8),
then we consider (1.5) at fixed $k$ as an equation for $H(k,\cdot)\in
L_{\mu}^{\infty}(\R^3)$. An analysis of equation (1.5)
 for $d=3$ and with (1.8) taken as a basic assumption on $v$ is given in
Section 3.

Note that, actually, $h$ on $\Theta$ is a zero energy restriction of a
function $h$ introduced by Faddeev (see [F2], [HN]) as an extention to the
complex domain of the classical scattering amplitude for the Schr\"odinger
equation at positive energies. In addition, the restriction
$h\big|_{\Theta}$ was not considered in Faddeev's works. Note that
$h\big|_{\Theta}$  was considered for the first time in [BC1] for $d=3$
in the framework of Problem 1a formulated below. The Faddeev function
$h$ was, actually, rediscovered in [BC1].
The fact that $\bar\pa$- scattering data of [BC1] coincide with the
Faddeev function $h$ was observed, in particular, in [HN].

In the present work, in addition to $h$ on $\Theta$, we consider
$h\big|_{\Gamma}$, $h\big|_{\Theta^{\tau}}$ and  $h\big|_{\Gamma^{\tau}}$,
where
$$\Gamma=\{k={p\over 2}+{i|p|\over 2}\gamma(p),\ \
l=-{p\over 2}+{i|p|\over 2}\gamma(p):\ p\in\R^d\},\eqno(1.10a)$$
where $\gamma$ is a piecewise continuous (or just measurable) function of
$p\in\R^d$ with values in $\S^{d-1}$ and such that
$$\gamma(p)p=0,\ \ p\in\R^d,\eqno(1.10b)$$
$$\eqalignno{
&\Theta^{\tau}=\{(k,l)\in\Theta:\ |Im\,k|=|Im\,l|<\tau\},&(1.11)\cr
&\Gamma^{\tau}=\Gamma\cap\Theta^{\tau},&(1.12)\cr}$$
where $\tau>0$. Note that
$$\eqalignno{
&\Gamma\subset\Theta,&(1.13)\cr
&dim\,\Theta=3d-4,\ \ dim\,\Gamma=dim\,\R^d=d,&(1.14)\cr
&3d-4=d\ \ {\rm for}\ \ d=2,\ \ 3d-4>d\ \ {\rm for}\ \ d\ge 3.&(1.15)\cr}$$

Using (1.4), (1.5) one can see that
$$h(k,l)\approx\hat v(p),\ \ (k,l)\in\Theta,\ \ k-l=p,\eqno(1.16)$$
in the Born approximation (that is in the linear approximation near zero
potential). Using (1.10), (1.13), (1.14), (1.16) one can see that, in
general, $h\big|_{\Gamma}$ is a nonlinear analog of the Fourier
transform $\hat v$. Note also that  $h\big|_{\Gamma}$  is a zero energy
analog of the reflection coefficient (backscattering amplitude) considered
(in particular) in [Mos], [P], [HN], [ER].

In the present work we consider, in particular, the following inverse
scattering problems for equation (1.1) under assumptions (1.2).

\vskip 2 mm
{\bf Problem 1.} (a)\ Given $h$ on $\Theta$, find $v$ on $\R^d$ (and
characterize $h$ on $\Theta$);

(b)\ Given $h$ on $\Theta^{\tau}$ for some (sufficiently great)
$\tau>0$, find $v$ on $\R^d$, at least, approximately.

\vskip 2 mm
{\bf Problem 2.} (a)\ Given $h$ on $\Gamma$, find $v$ on $\R^d$ (and
characterize $h$ on $\Gamma$);

(b)\ Given $h$ on $\Gamma^{\tau}$ for some (sufficiently great)
$\tau>0$, find $v$ on $\R^d$, at least, approximately.

Using (1.14), (1.15), (1.16) one can see that Problems 1a,1b are strongly
overdetermined for $d\ge 3$, whereas Problems 2a, 2b are nonoverdetermined
for $d\ge 2$ (at least, in the sense of the dimension considerations and
in the Born approximation). In addition, using (1.12), (1.13) one can see
that any reconstruction method for Problems 2 is also a reconstruction
method for Problems 1. The present work is focused on Problems 2a, 2b for
the most important three-dimensional case $d=3$.
In addition, we are focused on potentials $v$ with
$$\eqalign{
&\hat v\in L_{\mu}^{\infty}(\R^3)\ \ {\rm with\ sufficiently\ small}\ \
\|\hat v\|_{\mu}\cr
&{\rm for\ some\ fixed}\ \ \mu\ge 2, \cr}\eqno(1.17)$$
where $L_{\mu}^{\infty}(\R^3)$ and $\|\cdot\|_{\mu}$ are defined in (1.9).
In some results we also still assume for simplicity that
$\hat v\in {\cal C}(\R^3)$ (in addition to (1.8) or (1.17)), where
$\cal C$ denotes the space of continuous functions.
The main results of the present work include, in
particular:
\item{(  I)} uniqueness theorem, reconstruction procedure and stability
estimate for Problem 2a for $v$ satisfying (1.17)  (with $\hat v\in
{\cal C}(\R^3)$) (see Theorem 2.1) and
\item{( II)} approximate reconstruction method for Problem 2b for $v$
satisfying (1.17)  (with $\hat v\in {\cal C}(\R^3)$)
(see Theorem 2.1 and Corollary 2.1).
These results are formulated and proved in Sections 2-12. In the present
work we formulate also:
\item{(III)} characterization for Problem 2a for $v$ satisfying (1.17)
(see Theorem 2.2) and
\item{(IV)} new characterization for Problem 1a or more precisely a
characterization for Problem 1a for $v$ satisfying (1.17) (see Theorem 2.3).

We plan to give a complete proof of these characterizations in a separate
work, where we plan to show also that the aforementioned results I and II
remain
valid without the additional assumption that $\hat v\in {\cal C}(\R^3)$. All
these results I, II, III and IV are presented in detail in Section 2.

Note that Problem 1a was considered for the first time in [BC1] for $d=3$
from pure mathematical point of view without any physical applications.
No possibility to measure $h$ on $\Theta\b\{(0,0)\}$ directly in some
physical experiment is known at present. However, as it was shown in [No1]
(see also [HN] (Note added in proof), [Na1], [Na2], [No4]), Problems 1
naturally arise in the electrical impedance tomography and, more
generally, in the inverse boundary value problem (Problem 3) formulated as
follows. Consider the equation
$$-\Delta\psi+v(x)\psi=0,\ \ x\in D,\eqno(1.18)$$
where
$$\eqalign{
&D\ \ {\rm is\ an\ open\ bounded\ domain\ in}\ \ \R^d,\ \ d\ge 2,\cr
&{\rm with\ sufficiently\ regular\ boundary}\ \ \pa D,\cr
&v\ \ {\rm is\ a\ sufficiently\ regular\ function\ on}\ \ \bar D=D\cup\pa D.
\cr}\eqno(1.19)$$
We assume also that
$$\eqalign{
&0\ \ {\rm is\ not\ a\ Dirichlet\ eigenvalue\ for} \cr
&{\rm the\ operator}\ \ -\Delta+v\ \ {\rm in}\ \ D.\cr}\eqno(1.20)$$
Consider the map $\Phi$ such that
$${\pa\psi\over \pa\nu}\big|_{\pa D}=\Phi\bigl(\psi\big|_{\pa D}\bigr)
\eqno(1.21)$$
for all sufficiently regular solutions of (1.18) in $\bar D$, where $\nu$
is the outward normal to $\pa D$. The map $\Phi$ is called the
Dirichlet-to-Neumann map for equation (1.18). The aforementioned inverse
boundary value problem is:

\vskip 2 mm
{\bf Problem 3.}
Given $\Phi$, find $v$ on $D$.

In addition, the simplest  interpretation of $D$, $v$ and $\Phi$ in the
framework of the electrical impedance tomography consists in the
following (see [SU], [No1], [Na1]):\ $D$ is a body with isotropic
conductivity $\sigma(x)$ (where $\sigma\ge\sigma_{min}>0$),
$$\eqalignno{
&v(x)=(\sigma(x))^{-1/2}\,\Delta\,(\sigma(x))^{1/2},\ \ x\in D,&(1.22)\cr
&\Phi=\sigma^{-1/2}\bigl(\Lambda\sigma^{-1/2}+{\pa\sigma^{1/2}\over \pa\nu}
\bigr),&(1.23)\cr}$$
where $\Lambda$ is the voltage-to-current map on $\pa D$ and $\sigma^{1/2}$,
$\pa\sigma^{1/2}\big/\pa\nu$ in (1.23) denote the multiplication
operators by the functions $\sigma^{-1/2}\big|_{\pa D}$,
$\bigl(\pa\sigma^{1/2}\big/\pa\nu\bigr)\big|_{\pa D}$, respectively.

Note that the formulation of Problem 3 goes back to Gelfand [G] and
Calderon [C].

Returning to Problems 1, 2 and their relation to Problem 3 one can see that
the Faddeev function $h$ of Problems 1, 2 does not appear in Problem 3.
However, as it was shown in [No1] (see also [HN] (where this result of
[No1] was announced in Note added in proof),
[Na1], [Na2], [No4]), if $h$ corresponds to equation (1.1), where
$$\eqalign{
&v\ \ {\rm of}\ \ (1.1)\ \ {\rm coincides\ on}\ \ D\ \ {\rm with}\ \ v\ \
{\rm of}\ \ (1.18)\cr
&{\rm and}\ \ v\ \ {\rm of}\ \ (1.1)\ \ {\rm is\ identically\ zero\ on}\ \
\R^d\b\bar D,\cr}\eqno(1.24)$$
then $h$ on $\Theta$ can be determined from the Dirichlet-to-Neumann map
$\Phi$ for equation (1.18) via the following formulas and equation:
$$\eqalignno{
&h(k,l)=(2\pi)^{-d}\int\limits_{\pa D}\int\limits_{\pa D}e^{-ilx}
(\Phi-\Phi_0)(x,y)\psi(y,k)dydx\ \ {\rm for}\ \ (k,l)\in\Theta,&(1.25)\cr
&\psi(x,k)=e^{ikx}+\int\limits_{\pa D}A(x,y,k)\psi(y,k)dy,\ \ x\in\pa D,
&(1.26)\cr
&A(x,y,k)=\int\limits_{\pa D}G(x-z,k)(\Phi-\Phi_0)(z,y)dz,\ \ x,y\in\pa D,
&(1.27)\cr
&G(x,k)=-(2\pi)^{-d}e^{ikx}\int\limits_{\R^d}{e^{i\xi x}d\xi\over
{\xi^2+2k\xi}},\ \ x\in\R^d,&(1.28)\cr}$$
where $k\in\C^d$, $k^2=0$ in (1.26)-(1.28), $\Phi_0$ denotes the
Dirichlet-to-Neumann map for equation (1.18) for $v\equiv 0$, and
$(\Phi-\Phi_0)(x,y)$ is the Schwartz kernel of the integral operator
$\Phi-\Phi_0$. Note that (1.25), (1.27), (1.28) are explicit formulas,
whereas (1.26) is a linear integral equation (with parameter $k$) for
$\psi$ on $\pa D$. In addition, $G$ of (1.28) is the Faddeev's Green
function of [F1] for the Laplacian $\Delta$. Note also that formulas and
equation (1.25)- (1.27) are obtained and analyzed in [No1] for (1.19)
specified as
$$\eqalign{
&D\ \ {\rm is\ an\ open\ bounded\ domain\ in}\ \ \R^d,\ \ d\ge 2,\cr
&\pa D\in C^2,\ \ v\in L^{\infty}(D).\cr}\eqno(1.29)$$
Formulas and equation (1.25)-(1.27) reduce Problem 3 to Problems
1, 2. In addition, from numerical point of view $h(k,l)$ for
$(k,l)\in\Theta^{\tau}$ can be relatively easily determined from $\Phi$
via (1.27), (1.26), (1.25) if $\tau$ is sufficiently small. However,
if $(k,l)\in\Theta\b\Theta^{\tau}$, where $\tau$ is sufficiently great,
then the  determination of $h(k,l)$ from $\Phi$ via (1.27), (1.26),
(1.25) is very unstable (especially on the step (1.26)). The reason of
this instability is that formulas and equation (1.25)-(1.28) involve
the exponential functions $e^{-ilx}$, $e^{ikx}$ and, actually,
$e^{ik(x-z)}$ (arising in (1.27) in view of (1.28)), where
$(k,l)\in\Theta$, $x\in\pa D$, $z\in\pa D$, which  rapidly oscillate
in $x,z$ and may have exponentially great absolute values if
$(k,l)\in\Theta\b\Theta^{\tau}$ (and, therefore, $|Re\,k|=|Im\,k|=
|Re\,l|=|Im\,l|>\tau$) for sufficiently great $\tau$.

These remarks show that Problems 1, 2 are especially important in their
versions 1b, 2b as regards their applications to Problem 3 via (1.25)-(1.28)
(or via similar reductions). In addition, in view of (1.13)-(1.15), one
can see that it is much simpler to determine $h$ on $\Gamma$ (or on
$\Gamma^{\tau}$) only than completely on $\Theta$ (on on $\Theta^{\tau}$,
respectively) from $\Phi$ via (1.25)-(1.28) for $d\ge 3$. Therefore,
Problem 2b is of particular interest and importance in the framework of
applications of Problems 1, 2 to Problem 3 for $d\ge 3$.

In the present work we consider, mainly, Problems 1 and 2 for $d=3$. In
addition, as it was already mentioned, we are focused on
nonoverdetermined Problems 2a, 2b for $v$ satisfying (1.17). The main results
of the present work are presented in Section 2. (Some of these results
were already mentioned above.) Note that only restrictions in time prevent us
from generalizing all main results of the present work to the case $d>3$.
Actually, the results of the present work are obtained in the framework of a
development of the $\bar\pa$-approach to inverse scattering at fixed
energy in dimension $d\ge 3$ of [BC1], [HN], [No3], [No5].
In particular, the central part of the present work consists in an analysis
of the non-linear $\bar\pa$-equation (3.13) for the Faddeev function $H$ on
$\Theta$ for $v$ satisfying (1.17) (with $\hat v\in {\cal C}(\R^3)$), see
Sections 5,6,7.

Actually, in the present work we do not consider Problems 1 and 2 for $d=2$:
inverse scattering at fixed energy in dimension $d=2$ differs considerably
from inverse scattering at fixed energy in dimension $d\ge 3$. Note that a
global reconstruction method for Problem 2a for $d=2$ and for $v$ of the
form (1.2), (1.22), where $x\in\R^2$, $\sigma\ge\sigma_{min}>0$, was given
in [Na2] in the framework of a development of the $\bar\pa$-approach to
inverse scattering at fixed energy in dimension $d=2$ (see references to
[BLMP], [GN], [No2], [T] given in [Na2] in connection with this approach).
In addition, this result on Problem 2a is given in [Na2] in the
framework of applications to the (two-dimensional) electrical impedance
tomography via the reduction (1.25)- (1.27) for $d=2$ (given first in
[No1]). Besides, note that there is an essential similarity between the
results of [Na2] on global reconstruction for Problem 2a for $d=2$ and for
$v$ of the form (1.2), (1.22), where $x\in\R^2$, $\sigma\ge\sigma_{min}>0$,
and results of [BC2] on global inverse scattering reconstruction for some
$2\times 2$ first order system on the plane (see also [BU] in this
connection).

Applications of result of the present work to the electrical impedance
tomography and more generally to Problem 3 will be analyzed in detail in a
subsequent paper (where we plan to give, in particular, new stability
estimates for Problem 3).

Concerning results given in the literature on Problem 3, see [KV],
[SU], [HN]
(note added in proof), [No1], [A], [Na1], [Na2], [BU], [Ma], [No4]
and references therein.

\vskip 2 mm
{\bf 2. Main results}

As it was already mentioned in the introduction, the main results of the
present work include, in particular:
\item{( I)} uniqueness theorem, reconstruction procedure and stability
estimate for Problem 2a for $v$ satisfying (1.17) (with $\hat v\in
{\cal C}(\R^3)$) and
\item{(II)} approximate reconstruction method for Problem 2b for $v$
satisfying (1.17)  (with $\hat v\in {\cal C}(\R^3)$),

\noindent
see Theorem 2.1 and
Corollary 2.1 formulated below in this section (and proved by means of
analysis developed in Sections 3-12).

We identify $h\big|_{\Gamma}$ and $h\big|_{\Gamma^{\tau}}$ with $R$ and
$R_{2\tau}$ on $\R^d$, where
$$R(p)=h\bigl({p\over 2}+{i|p|\over 2}\gamma(p),-
{p\over 2}+{i|p|\over 2}\gamma(p)\bigr),\ \ p\in\R^d,\eqno(2.1)$$
$$\eqalign{
&R_{2\tau}(p)=R(p)\ \ {\rm for}\ \ |p|< 2\tau,\ \ p\in\R^d,\cr
&R_{2\tau}(p)=0\ \ {\rm for}\ \ |p|\ge 2\tau,\ \ p\in\R^d,\cr}\eqno(2.2)$$
where $\gamma$ is the function of (1.10).

\vskip 2 mm
{\bf Theorem 2.1.}
{\it Let}
$$\eqalignno{
&\hat v\in L_{\mu}^{\infty}(\R^3)\ \ {\it for\ some}\ \ \mu\ge 2,&(2.3)\cr
&\|\hat v\|_{\mu}\le C < {1\over {c_1(\mu)+8c_6(\mu)}},&(2.4)\cr}$$
{\it where} $L_{\mu}^{\infty}(\R^3)$ {\it and} $\|\cdot\|_{\mu}$ {\it are
defined in} (1.9), $c_1(\mu)$ {\it and} $c_6(\mu)$ {\it are the positive
constants of Lemmas 3.1 and} 6.1. {\it (For simplicity we also still assume
that} $\hat v\in {\cal C}(\R^3)$.)
{\it Let} $R$ {\it be defined by} (2.1)
({\it for
some given} $\gamma$ {\it of} (1.10) {\it for} $d=3$). {\it Then}
$$R\in L_{\mu}^{\infty}(\R^3),\ \ \|R\|_{\mu}\le {C\over {1-c_1(\mu)C}},
\eqno(2.5)$$
{\it and} $R$ {\it uniquely determines} $\hat v$ {\it via the following
reconstruction procedure}
$$R\buildrel (6.6) \over \rightarrow H\buildrel (7.1b),(7.2b) \over
\rightarrow \hat v,\eqno(2.6)$$
{\it where} (6.6) {\it is a nonlinear integral equation of Proposition} 6.1
{\it of Section} 6, (7.1b), (7.2b) {\it are explicit formulas of Section} 7
{\it and where we solve} (6.6) {\it by the method of successive
approximations (see Proposition} 6.2 {\it and Lemma} 6.2). {\it In
addition, if\ }\ $R_{appr}$ {\it is an arbitrary approximation to} $R$,
{\it where}  $R_{appr}$ {\it also satisfies} (2.5), {\it and}
${\hat v}_{appr}$ {\it is determined from} $R_{appr}$ {\it via} (2.6)
({\it with} $H$ {\it replaced by}  $H_{appr}$), {\it then the following
stability estimate holds}:
$$\|\hat v-\hat v_{appr}\|_{\mu}\le
{{1-c_1(\mu)C}\over {1-(c_1(\mu)+8c_6(\mu))C}} \|R-R_{appr}\|_{\mu}.
\eqno(2.7)$$

One can see that Theorem 2.1 includes uniqueness theorem, reconstruction
procedure and stability estimate for Problem 2a (of the introduction) for
$v$ satisfying (1.17) (with $\hat v\in {\cal C}(\R^3)$).

Theorem 2.1 follows from Proposition 3.1, Lemmas 6.2, 6.3, Propositions 6.1,
6.2 and
formulas (7.1), (7.2) (of Sections 3,6 and 7). In particular, condition (2.4)
of Theorem 2.1 implies condition (6.20) of Proposition 6.2 and condition
(3.6) of (part I of) Proposition 3.1.

\vskip 2 mm
{\bf Corollary 2.1.}
{\it Let} $v$ {\it satisfy} (2.3), (2.4) {\it and, in addition},
$$\hat v\in L_{\mu^*}^{\infty}(\R^3)\ \ {\it for\ some}\ \ \mu^*>\mu.
\eqno(2.8)$$
{\it (For simplicity we also still assume that} $\hat v\in {\cal C}(\R^3)$.)
{\it Let} $\hat v_{2\tau}$ {\it denotes} $\hat v_{appr}$ {\it
reconstructed from} $R_{appr}$ {\it via} (2.6) ({\it as in Theorem} 2.1),
{\it where} $R_{appr}=R_{2\tau}$ ({\it defined by} (2.1), (2.2) {\it for}
$d=3$). {\it Then}
$$R\in L_{\mu^*}^{\infty}(\R^3) \eqno(2.9)$$
{\it and}
$$\|\hat v-\hat v_{2\tau}\|_{\mu}\le
{{1-c_1(\mu)C}\over {1-(c_1(\mu)+8c_6(\mu))C}}
{\|R\|_{\mu^*}\over (1+2\tau)^{\mu^*-\mu}}\ \ {\it for}\ \ \tau>0.
\eqno(2.10)$$
One can see that Theorem 2.1 and Corollary 2.1 give an approximate
reconstruction method for Problem 2b (of the introduction) for $v$
satisfying (1.17) (with $\hat v\in {\cal C}(\R^3)$).

\vskip 2 mm
Note that (2.9) follows from the property that $R\in L_{\mu}^{\infty}(\R^3)$,
the assumption (2.8) and the part II of Proposition 3.1 with
$\mu=\mu^*$. Further,
Corollary 2.1 follows from Theorem 2.1 and estimates (6.25), (6.26).
The approximate reconstruction of Corollary 2.1 is presented in more
detail in Proposition 6.3 complemented by formulas (7.5)-(7.8).

One can see that Theorem 2.1 and Corollary 2.1 give also reconstruction
results for Problem 1a and Problem 1b ( of the  introduction) for $d=3$
and $v$ satisfying (1.17) (with $\hat v\in {\cal C}(\R^3)$).
Let us compare these results with the
reconstructions for Problems 1a and 1b for $d=3$ via formulas (2.11),
(2.12), (2.14) presented below. From formula (1.4), equation (1.5) and
Proposition 3.1  (of Section 3) it follows
that if $v$ satisfies (2.3), then
$$\hat v(p)=\lim_{\scriptstyle (k,l)\in\Theta,\ k-l=p \atop\scriptstyle
|Im\,k|=|Im\,l|=\tau\to\infty}h(k,l)\ \ {\rm for\ any}\ \ p\in\R^3,
\eqno(2.11)$$
$$\eqalign{
&|\hat v(p)-h(k,l)|\le {2c_2(\mu)C^2\over (1+|p|)^{\mu}}\,
{(\ln\,\tau)^2\over \tau}\cr
&{\rm for}\ \ (k,l)\in\Theta,\ \ p=k-l,\ \
|Im\,k|=|Im\,l|=\tau\ge \tau(C,\mu),\ \ \|\hat v\|_{\mu}\le C,\cr}
\eqno(2.12)$$
where $c_2(\mu)$ is the constant of Lemma 3.1 and $\tau(C,\mu)$ is the
smallest number such that
$$c_2(\mu)C{(\ln\,\tau(C,\mu))^2\over \tau(C,\mu)}\le {1\over 2},\ \
\ln\,\tau(C,\mu)\ge 2.$$
Actually, for sufficiently regular $v$ on $\R^3$ with sufficient decay at
infinity formula (2.11) and some results of the type (2.12) (with less
precise right-hand side) were given first in [HN]. Note also that if
$$\eqalign{
&v\in L^{\infty}(\R^3),\ \ ess\,\sup_{x\in\R^3}(1+|x|)^{3+\ep}|v(x)|\le C \cr
&{\rm for\ some\ positive}\ \ \ep\ \ {\rm and}\ \ C,\cr}\eqno(2.13)$$
then
$$\eqalign{
&|\hat v(p)-h(k,l)|\le {2\tilde c_2(\ep) C^2\over \tau}\ \ {\rm for}\ \
(k,l)\in\Theta,\ \ p=k-l,\cr
&|Im\,k|=|Im\,l|=\tau\ge\tilde\tau(C,\ep),\cr}\eqno(2.14)$$
where $\tilde c_2(\ep)$ and $\tilde\tau(C,\ep)$ are some positive constants
(similar to constants $c_2(\mu)$ and $\tau(C,\mu)$ of (2.12)) (see [Na1] and
[No3] as regards estimate  (2.14) under assumption (2.13)). One can see that
for $d=3$ already the simple formulas (2.11), (2.12), (2.14) give a
reconstruction method for Problem 1a and an approximate reconstruction
method for Problem 1b. However, for this approximate reconstruction of
the Fourier transform $\hat v$ from $h$ on $\Theta^{\tau}$ via (2.12),
(2.14) the error decaies rather slowly as $\tau\to +\infty$: even for $v$
of the Schwartz class on $\R^3$ the decay rate of this error, for example,
in the uniform norm on the ball ${\cal B}_{\tau}=\{p\in\R^3:\ |p|\le r\}$,
where $r>0$ is fixed, is not faster than $O(1/\tau)$ as $\tau\to +\infty$.
An important advantage of the approximation $\hat v_{2\tau}$ of
Corollary 2.1 in comparison with the approximate reconstruction based on
(2.12), (2.14) consists in a fast decay of the error norm
$\|\hat v-\hat v_{2\tau}\|_{\mu}=O(1/\tau^{\mu^*-\mu})$ as $\tau\to +\infty$
(see estimate (2.10)), at least, if $\mu^*-\mu$ is sufficiently great. For
example, if $v$ belongs to the Schwartz class on $\R^3$ and, as in
Theorem 2.1 and Corollary 2.1, is sufficiently small in the  sense (2.4)
for some $\mu$, then estimate (2.10) holds for any $\mu^*>\mu$ and
$\|\hat v-\hat v_{2\tau}\|_{\mu}=O(\tau^{-\infty})$ as $\tau\to +\infty$.
This fast convergence of $\hat v_{2\tau}$ to $\hat v$ as $\tau\to +\infty$
is in particular important in the framework of applications to Problem 3
(of the introduction) via the reduction (1.25)-(1.27): the point is that
the determination of $h\big|_{\Theta^{\tau}}$ from $\Phi$ via
(1.25)-(1.27) is sufficiently stable for sufficiently small $\tau$ only
(see related discussion of the introduction), but  $\hat v_{2\tau}$
reconstructed from $h\big|_{\Gamma^{\tau}}$ (as described in Corollary 2.1)
 well approximates $\hat v$ even if $\tau$  is relatively small (due to
the rapid decay of the error
$\hat v-\hat v_{2\tau}$ as $\tau\to +\infty$).
An obvious disadvantage of Theorem 2.1 and Corollary 2.1 in comparison with
formulas (2.11), (2.12), (2.14) consists in the small norm assumption (2.4).
In a subsequent work we plan to propose an approximate reconstruction of
$\hat v$ from $h$ on $\Theta^{\tau}$ (for $d=3$) with a similar (fast)
decay of the error for $\tau\to +\infty$ as in Corollary 2.1 but without the
assumption that $v$ is small in some sense.

As it was already mentioned in the introduction, in the present work we
formulate also:
\item{(III)} characterization for Problem 2a for $v$ satisfying (1.17) and
\item{( IV)} new characterization for Problem 1a or more precisely a
characterization for

\item{     } Problem 1a for  $v$ satisfying (1.17),

\noindent
see Theorems 2.2 and 2.3 presented next.

\vskip 2 mm
{\bf Theorem 2.2}.
{\it Let} $v$ {\it satisfy} (2.3) {\it and}
$$\|\hat v\|_{\mu}\le C<1/c_1(\mu),\eqno(2.15)$$
{\it where} $c_1(\mu)$ {\it is the constant of Lemma} 3.1. {\it Then} $R$
({\it defined according to} (2.1), (1.4), (1.5))
{\it satisfies} (2.5). {\it Conversely, let }
$$R\in L^{\infty}_{\mu}(\R^3)\ \ {\it for\ some}\ \ \mu\ge 2 \eqno(2.16)$$
{\it and}
$$\|R\|_{\mu}\le r/2,\ \ r< c_7(\mu),\eqno(2.17)$$
{\it where} $c_7(\mu)$ {\it is some positive constant.}
{\it Then} $R$ {\it is the scattering data (defined according to} (2.1),
(1.4), (1.5)) {\it for some potential} $v$, {\it where}
$$\hat v\in L^{\infty}_{\mu}(\R^3),\ \  \|\hat v\|_{\mu}\le r.\eqno(2.18)$$

One can see that Theorem 2.2 gives a characterization for Problem 2a (of the
introduction) for $v$ satisfying (1.17).

Consider
$$\eqalignno{
&\Omega=\{k\in\C^3,\ p\in\R^3:\ k^2=0,\ p^2=2kp\},&(2.19)\cr
&\Xi=\{(k,p):\ k={p\over 2}+{i|p|\over 2}\gamma(p),\ p\in\R^3\},&(2.20)\cr}
$$where $\gamma$ is the function of (1.10).

Note that
$$\Omega\approx\Theta,\ \ \Xi\approx\Gamma \eqno(2.21)$$
or more precisely
$$\eqalign{
&(k,p)\in\Omega\Rightarrow (k,k-p)\in\Theta,\ \
(k,l)\in\Theta\Rightarrow (k,k-l)\in\Omega,\cr
&(k,p)\in\Xi\Rightarrow (k,k-p)\in\Gamma,\ \
(k,l)\in\Gamma\Rightarrow (k,k-l)\in\Xi,\cr}\eqno(2.22)$$
where $\Theta$ and $\Gamma$ are defined by (1.3) and (1.10a) for $d=3$.
Due to (2.21), (2.22), $h$ on $\Theta$ in Problem 1 for $d=3$ can be
considered as
$H$ on $\Omega$  and $h$ on $\Gamma$ in Problem 2 for $d=3$ can be considered
as $H$ on $\Xi$, where $h$ and $H$ are related by (1.4).

Consider
$$\eqalign{
&L^{\infty}_{\mu}(\Omega)=\{U\in L^{\infty}(\Omega):\ |||U|||_{\mu}<+\infty
\},\cr
&|||U|||_{\mu}=ess\,\sup_{(k,p)\in\Omega}(1+|p|)^{\mu}|U(k,p)|,\ \
\mu>0.\cr}\eqno(2.23)$$

\vskip 2 mm
{\bf Theorem 2.3.}
{\it Let} $v$ {\it satisfy} (2.3), (2.15) {\it and} $H$ {\it be defined on}
$\Omega$ {\it by means of} (1.5). {\it Then}
$$H\in L^{\infty}_{\mu}(\Omega),\ \ |||H|||_{\mu}\le {C\over {1-c_1(\mu)C}},
\eqno(2.24)$$
{\it and for almost any} $p\in\R^3\b 0$ {\it the} $\bar\pa$-
{\it equation} (3.13) {\it for} $H$ {\it on} $\Omega$ {\it holds.}

{\it Conversely, let}
$$\eqalignno{
&H\in L^{\infty}_{\mu}(\Omega)\ \ {\it for\ some}\ \ \mu\ge 2,&(2.25)\cr
&|||H|||_{\mu}\le r,\ \ r< c_8(\mu),&(2.26)\cr}$$
{\it where} $c_8$ {\it is a positive constant, and for almost any}
$p\in\R^3\b 0$ {\it the} $\bar\pa$- {\it equation} (3.13) {\it holds. Then}
$H$ {\it on} $\Omega$ {\it is the scattering data} ({\it defined using}
(1.5)) {\it for some potential} $v$, {\it where}
$$\hat v\in L^{\infty}_{\mu}(\R^3),\ \ \|\hat v\|_{\mu}\le r.\eqno(2.27)$$
One can see that Theorem 2.3 gives a characterization for Problem 1a (of the
introduction) for $v$ satisfying (1.17) (and where $h$ on $\Theta$ is
considered as $H$ on $\Omega$).
In a separate work we plan to give a detailed comparison of Theorem 2.3 with
related results of [BC1] and [HN]. In particular, Theorem 2.3 develops
and simplifies the results of [BC] on the range characterization of $H$ on
$\Omega$.

The scheme of proof of Theorems 2.2 and 2.3 consists in the following:
\item{(1)} The result that (2.3), (2.15) imply (2.5) and (2.24) follows
from Proposition 3.1.
\item{(2)} It is a separate lemma that the $\bar\pa$- equation (3.13)
remains valid for almost any $p\in\R^3\b 0$ if $v$ satisfies (2.3) and
(2.15).
\item{(3)} To prove the sufficiency parts of Theorems 2.2 and 2.3, we use
Proposition 3.1, the aforementioned separate lemma concerning the
$\bar\pa$- equation (3.13), and the analysis developed in Sections 4,5,6 and
7. In addition, in the framework of this proof we obtain that the
constants $c_7(\mu)$ and $c_8(\mu)$  of Theorems 2.2 and 2.3 can be defined
as follows:
$$\eqalignno{
&c_7(\mu)=\min\,\bigl({1\over 4c_6(\mu)},\ \ {1\over {c_1(\mu)+2c_6(\mu)}}
\bigr),&(2.28)\cr
&c_8(\mu)={1\over {c_1(\mu)+2c_6(\mu)}},&(2.29)\cr}$$
where $c_1$ and $c_6$ are the constants of Lemmas 3.1 and 6.1.

On the basis of this scheme we plan to give a complete proof of Theorems 2.2
and 2.3 in a separate work, where we plan to show also that Theorem 2.1
and Corollary 2.1 remain valid without the additional assumption that
$\hat v\in {\cal C}(\R^3)$.

\vskip 2 mm
{\bf 3. Some results on direct scattering}
\vskip 2 mm
In this section we give some results on direct scattering at zero energy in
three dimensions or, more precisely, some results concerning equation (1.5)
and the function $H$ of (1.5) under assumption (1.8).

Consider the operator $A(k)$ from (1.5) for $d=3$:
$$(A(k)U)(p)=\int\limits_{\R^3}{\hat v(p+\xi)U(-\xi)d\xi\over {\xi^2+2k\xi}},
\ \ p\in\R^3,\ \ k\in\Sigma,\eqno(3.1)$$
where $U$ is a test function, $\Sigma$ is defined by (1.7) for $d=3$.
Let $\cal C$ stand for continuous functions.

\vskip 2 mm
{\bf Lemma 3.1.}
{\it Let} $v$ {\it satisfy} (1.8), $A(k)$ {\it be defined by} (3.1) {\it and}
$U\in L^{\infty}_{\mu}(\R^3)$. {\it Then:}
$$\eqalignno{
&A(k)U\in {\cal C}(\R^3),&(3.2)\cr
&\|A(k)U\|_{\mu}\le c_1(\mu)\|\hat v\|_{\mu}\|U\|_{\mu},&(3.3a)\cr
&\|A(k)U\|_{\mu}\le c_2(\mu)\|\hat v\|_{\mu}\|U\|_{\mu}
{(\ln\,(|Im\,k|))^2\over |Im\,k|},\ \ \ln\,|Im\,k|\ge 2,&(3.3b)\cr}$$
{\it for} $k\in\Sigma$ ({\it defined by} (1.7) {\it for}
$d=3$), {\it where} $c_1(\mu)$, $c_2(\mu)$ {\it and} $\rho(\mu)$ {\it are
some positive constants;  in addition},
$$\|(A(k^{\prime})-A(k))U\|_{\mu}\le\Delta(k,k^{\prime})\|\hat v\|_{\mu}
\|U\|_{\mu} \eqno(3.4a)$$
{\it for some} $\Delta(k,k^{\prime})$ {\it such that}
$$\lim_{k^{\prime}\to k}\Delta(k,k^{\prime})=0,\eqno(3.4b)$$
{\it where} $k,k^{\prime}\in\Sigma$; {\it in addition},
$$(A(k)U)(p)\in {\cal C}(\Sigma\times\R^3)\ \  {\it as\ a\ function\ of}\ \
k\ \ {\it and}\ \ p.\eqno(3.5)$$
Lemma 3.1 is proved in Section 8.

\vskip 2 mm
{\bf Proposition 3.1.}
{\it Let} $v$ {\it satisfy} (1.8) and  $\|\hat v\|_{\mu}\le C$. {\it Then the
 following statements are valid}:

\noindent
(I) {\it if}
$$\eta_1(C)\buildrel \rm def \over =c_1(\mu)C<1,\eqno(3.6)$$
{\it then equation} (1.5) {\it is uniquely solvable for} $H(k,\cdot)\in
L^{\infty}_{\mu}(\R^3)$ {\it for any} $k\in\Sigma$ ({\it by the method of
successive approximations}) {\it and}
$$\eqalignno{
&\|H(k,\cdot)\|_{\mu}\le {C\over {1-c_1(\mu)C}},\ \ k\in\Sigma,&(3.7)\cr
&H-\hat v\in {\cal C}(\Sigma\times\R^3),&(3.8a)\cr
&|H(k,p)-\hat v(p)|\le {c_1(\mu)C^2\over (1-c_1(\mu)C)(1+|p|)^{\mu}},\ \
k\in\Sigma,\ \ p\in\R^3;&(3.8b)\cr}$$

\noindent
(II) {\it if}
$$\eta_2(C,\tau)\buildrel \rm def \over =c_2(\mu)C{(\ln\tau)^2\over \tau}<1,
\ \ \ln\tau\ge 2,
\eqno(3.9)$$
{\it then equation} (1.5) {\it is uniquely solvable} ({\it by the method of
successive approximations}) {\it for} $H(k,\cdot)\in L^{\infty}_{\mu}(\R^3)$
{\it for any} $k\in\Sigma\b\Sigma^{\tau}$, {\it where}
$$\Sigma^{\tau}=\{k\in\Sigma:\ |Im\,k|<\tau\},\eqno(3.10)$$
{\it and}
$$\eqalignno{
&|H(k,\cdot)|_{\mu}\le {C\over {1-\eta_2(C,|Im\,k|)}},\ \ k\in\Sigma\b
\Sigma^{\tau},&(3.11)\cr
&H-\hat v\in {\cal C}((\Sigma\b\Sigma^{\tau})\times\R^3),&(3.12a)\cr
&|H(k,p)-\hat v(p)|\le {\eta_2(C,|Im\,k|)C\over (1-\eta_2(C,|Im\,k|))
(1+|p|)^{\mu}},\ \ k\in\Sigma\b\Sigma^{\tau},\ \ p\in\R^3.&(3.12b)\cr}$$

Proposition 3.1 is proved in Section 8.

Further, note that if $v$ satisfies (1.8) and $\|\hat v\|_{\mu}\le C$, where
$C$ satisfies (3.6), and also $\hat v\in {\cal C}(\R^3)$, then the Faddeev
function $H$ (of the part I of Proposition 3.1) satisfies the following
$\bar\pa$- equation on $\Omega$:
$$\eqalign{
&\bar\pa_kH(k,p)\big|_{Z_p}=\cr
&\sum_{j=1}^3\biggl(-2\pi\int\limits_{\xi\in S_k}\xi_jH(k,-\xi)H(k+\xi,p+\xi)
{ds\over |Im\,k|^2}\biggr)d\bar k_j\big|_{Z_p}\cr}\eqno(3.13)$$
for any $p\in\R^3\b 0$,  where
$$\eqalignno{
&Z_p=\{k\in\C^3:\ (k,p)\in\Omega\},\ \ p\in\R^3\b 0,&(3.14)\cr
&S_k=\{\xi\in\R^3:\ \xi^2+2k\xi=0\},\ \ k\in Z_p,&(3.15)\cr}$$
$ds$ is arc-length measure on the circle $S_k$  in $\R^3$.
Note also that, under the assumptions of the part II of Proposition 3.1
with $\hat v\in {\cal C}(\R^3)$, the $\bar\pa$ -equation (3.13) remains
valid with $Z_p$ replaced by $Z_p\cap (\Sigma\b\Sigma^{\tau})$.
Actually, at least under somewhat stronger assumptions on $v$ than in the
part I of Proposition 3.1  with $\hat v\in {\cal C}(\R^3)$, the $\bar\pa$ -
equation (3.13) was obtained for the first time in [BC1].

\vskip 2 mm
{\bf 4. Coordinates on $\Omega$}
\vskip 2 mm
Consider $\Omega$ defined by (2.19). For our considerations we introduce
some convinient coordinates on $\Omega$. Let
$$\Omega_{\nu}=\{k\in\C^3,\ p\in\R^3\b {\cal L}_{\nu}:\ k^2=0,\ p^2=2kp\},
\eqno(4.1)$$
where
$${\cal L}_{\nu}=\{p\in\R^3:\ p=t\nu,\ t\in\R\},\ \ \nu\in\S^2.\eqno(4.2)$$
Note that $\Omega_{\nu}$ is an open and dense subset of $\Omega$.

For $p\in\R^3\b {\cal L}_{\nu}$ consider $\theta(p)$ and $\omega(p)$ such
that
$$\eqalign{
&\theta(p),\omega(p)\ \ {\rm smoothly\ depend\ on}\ \ p\in\R^3\b
{\cal L}_{\nu},\cr
&{\rm take\ values\ in}\ \ \S^2,\ \ {\rm and}\cr
&\theta(p)p=0,\ \omega(p)p=0,\ \theta(p)\omega(p)=0.\cr}\eqno(4.3)$$
Note that (4.3) implies that
$$\omega(p)={p\times\theta(p)\over |p|}\ \ {\rm for}\ \
p\in\R^3\b {\cal L}_{\nu} \eqno(4.4a)$$
or
$$\omega(p)=-{p\times\theta(p)\over |p|}\ \ {\rm for}\ \
p\in\R^3\b {\cal L}_{\nu}, \eqno(4.4b)$$
where $\times$ denotes vector product.

To satisfy (4.3), (4.4a) we can take
$$\theta(p)={\nu\times p\over |\nu\times p|},\ \omega(p)=
{p\times\theta(p)\over |p|},\ p\in\R^3\b {\cal L}_{\nu}.\eqno(4.5)$$

\vskip 2 mm
{\bf Lemma 4.1.}
{\it Let} $\theta,\omega$ {\it satisfy} (4.3). {\it Then the following
formulas give a diffeomorphism between} $\Omega_{\nu}$ {\it and}
$(\C\b 0)\times (\R^3\b {\cal L}_{\nu})$:
$$\eqalignno{
&(k,p)\to (\lambda,p),\ \ {\it where}\ \
\lambda=\lambda(k,p)={2k(\theta(p)+i\omega(p))\over i|p|},&(4.6)\cr
&(\lambda,p)\to (k,p),\ \ {\it where}\ \ k=k(\lambda,p)=\kappa_1(\lambda,p)
\theta(p)+\kappa_2(\lambda,p)\omega(p)+{p\over 2},&(4.7)\cr}$$
$$\kappa_1(\lambda,p)={i|p|\over 4}(\lambda+{1\over \lambda}),
\ \  \kappa_2(\lambda,p)={|p|\over 4}(\lambda-{1\over \lambda}),$$
{\it where} $(k,p)\in\Omega_{\nu}$,
$(\lambda,p)\in (\C\b 0)\times (\R^3\b {\cal L}_{\nu})$.

Actually, Lemma 4.1 follows from properties (4.3) and the result that
formulas (4.6),

\noindent
(4.7) for $\lambda(k)$ and $k(\lambda)$ at fixed
$p\in\R^3\b {\cal L}_{\nu}$ give a diffeomorphism between $\{k\in\C^3:\
k^2=0,\ p^2=2kp\}$ and $\C\b 0$. The latter result follows from the
fact (see [GN],[No2]) that the following formulas
$$\lambda={{k_1+ik_2}\over i|E|^{1/2}},
\ \ k_1={i|E|^{1/2}\over 2}\bigl(\lambda+{1\over \lambda}\bigr),\ \
k_2={|E|^{1/2}\over 2}\bigl(\lambda-{1\over \lambda}\bigr)$$
give a diffeomorphism between $\{k\in\C^2:\ k^2=E\}$, $E<0$, and $\C\b 0$.

Note that for $k$ and $\lambda$ of (4.6), (4.7) the following formulas
hold:
$$|Im\,k|={|p|\over 4}\bigl(|\lambda|+{1\over |\lambda|}\bigr),\ \
|Re\,k|={|p|\over 4}\bigl(|\lambda|+{1\over |\lambda|}\bigr),\eqno(4.8)$$
where $(k,p)\in\Omega_{\nu}$,
$(\lambda,p)\in (\C\b 0)\times (\R^3\b {\cal L}_{\nu})$.

We consider $\lambda,p$ of Lemma 4.1 as coordinates on $\Omega_{\nu}$ and
on $\Omega$.

\vskip 2 mm
{\bf 5. $\bar\pa$-equation for $H$ on $\Omega$ in the coordinates $\lambda$,
$p$}
\vskip 2 mm
{\bf Lemma 5.1.}
{\it Let the assumptions of the part I of Proposition} 3.1 {\it be fulfilled
and} $\hat v\in {\cal C}(\R^3)$. {\it Let} $\lambda$, $p$ {\it be the
coordinates of Lemma} 4.1, {\it where} $\theta$, $\omega$ {\it satisfy}
(4.3), (4.4a). {\it Then}
$$\eqalign{
&{\pa\over \pa\bar\lambda}H(k(\lambda,p),p)=-{\pi\over 4}
\int_{-\pi}^{\pi}\biggl({|p|\over 2} {(|\lambda|^2-1)\over \bar\lambda
|\lambda|}(\cos\v-1)-|p|{1\over \bar\lambda}\sin\v\biggr)\times\cr
&H(k(\lambda,p),-\xi(\lambda,p,\v))H(k(\lambda,p,\v),p+\xi(\lambda,p,\v))
d\v\cr}\eqno(5.1)$$
{\it for} $\lambda\in\C\b 0$,   $p\in\R^3\b {\cal L}_{\nu}$, {\it where}
$k(\lambda,p)$
{\it is defined in} (4.7) ({\it and also depends on} $\nu$, $\theta$,
$\omega$),
$$\eqalignno{
&\xi(\lambda,p,\v)=Re\,k(\lambda,p)(\cos\v-1)+k^{\perp}(\lambda,p)\sin\v,
&(5.2)\cr
&k^{\perp}(\lambda,p)={Im\,k(\lambda,p)\times Re\,k(\lambda,p)\over
|Im\,k(\lambda,p)|},&(5.3)\cr}$$
{\it where} $\times$ {\it in} (5.3) {\it denotes vector product}.

Proof of Lemma 5.1 is given in Section 9. In this proof we deduce (5.1)
from (3.13).

Note that (5.1) can be written as
$${\pa\over \pa\bar\lambda}H(k(\lambda,p),p)=\{H,H\}(\lambda,p),\ \
\lambda\in\C\b 0,\ \  p\in\R^3\b {\cal L}_{\nu},\eqno(5.4)$$
where
$$\eqalign{
&\{U_1,U_2\}(\lambda,p)=-{\pi\over 4}
\int_{-\pi}^{\pi}\biggl({|p|\over 2} {{|\lambda|^2-1}\over \bar\lambda
|\lambda|}(\cos\v-1)-{|p|\over \bar\lambda}\sin\v\biggr)\times\cr
&U_1(k(\lambda,p),-\xi(\lambda,p,\v))U_2(k(\lambda,p)+\xi(\lambda,p,\v),
p+\xi(\lambda,p,\v))d\v,\cr}\eqno(5.5)$$
where $U_1$, $U_2$ are test functions on $\Omega$ (defined by (2.19)) and
$k(\lambda,p)$, $\xi(\lambda,p,\v)$ are defined by (4.7), (5.2),
$(\lambda,p)\in (\C\b 0)\times (\R^3\b {\cal L}_{\nu})$.
Note that in the left-hand side of (5.1), (5.4)
$$(k(\lambda,p),p)\in\Omega_{\nu} \eqno(5.6a)$$
and in the right-hand side of (5.1), (5.5)
$$\eqalign{
&(k(\lambda,p),-\xi(\lambda,p,\v))\in\Omega\b (0,0),\cr
&(k(\lambda,p)+\xi(\lambda,p,\v),p+\xi(\lambda,p,\v))\in\Omega\b (0,0),\cr}$$
where $\lambda\in\C\b 0$, $p\in\R^3\b {\cal L}_{\nu}$, $\v\in [-\pi,\pi]$
(and (0,0) denotes the point $\{k=0,p=0\}$).

\vskip 2 mm
{\bf Lemma 5.2.}
{\it Let the assumptions of Lemma} 4.1 {\it be fulfilled. Let}
$U_1,U_2\in L^{\infty}_{\mu}(\Omega)$
{\it for some} $\mu\ge 2$, {\it where}  $L^{\infty}_{\mu}(\Omega)$ {\it is
defined by} (2.23). {\it Let} $\{U_1,U_2\}$ {\it be defined by}
(5.5). {\it Then}:
$$\{U_1,U_2\}\in L^{\infty}((\C\b 0)\times (\R^3\b {\cal L}_{\nu})) \eqno(5.7)
$$
{\it and}
$$\eqalign{
&|\{U_1,U_2\}(\lambda,p)|\le {|||U_1|||_{\mu}|||U_2|||_{\mu}\over
(1+|p|)^{\mu}}\times\cr
&\biggl({c_3(\mu)|\lambda|\over (|\lambda|^2+1)^2}+
{c_4(\mu)|p|||\lambda|^2-1|\over |\lambda|^2(1+|p|(|\lambda|+|\lambda|^{-1}))
^2}+
{c_5(\mu)|p|\over |\lambda|(1+|p|(|\lambda|+|\lambda|^{-1}))}\biggr)\cr}
\eqno(5.8)$$
{\it for almost all}
$(\lambda,p)\in (\C\b 0)\times (\R^3\b {\cal L}_{\nu})$.

Proof of Lemma 5.2 is given in Section 10.

\vskip 2 mm
{\bf 6. Finding $H$ on $\Omega$ from its nonredundant restrictions
$H\big|_{\Xi}$}
\vskip 2 mm
Our next purpose is to give an integral equation for finding $H$  on $\Omega$
from $R=H\big|_{\Xi}$, where $\Omega$ and $\Xi$ are defined by (2.19),
(2.20). Actually, we will give an integral equation for finding $H$ on
$\Omega_{\nu}$ from $R=H\big|_{\Xi_{\nu}}$, where $\Omega_{\nu}$ is
defined by (4.1) and $\Xi_{\nu}=\Xi\cap\Omega_{\nu}$. In the coordinates
of Lemma 4.1 this means that we will give an integral equation for
finding
$$H(\lambda,p)=H(k(\lambda,p),p),\ \ \lambda\in\C\b 0,\ \
p\in\R^3\b {\cal L}_{\nu},\eqno(6.1)$$
from
$$R(p)=H(\lambda_0(p),p)=H(k(\lambda_0(p),p),p),\ \
p\in\R^3\b {\cal L}_{\nu},\eqno(6.2)$$
where $\lambda_0$ of (6.2) is a  piecewise continuous
function of  $p\in\R^3\b {\cal L}_{\nu}$
with values in
$$T=\{\lambda\in\C:\ |\lambda|=1\}.\eqno(6.3)$$
These properties of $\lambda_0$ of (6.2) follow from the properties of
$\gamma$ of (1.10a) and from (4.6). Note that if, for example,
$\gamma=\theta$, where $\theta$, $\omega$ are defined by (4.5), then
$\lambda_0(p)\equiv 1$ for  $p\in\R^3\b {\cal L}_{\nu}$.

We will use the following formula
$$\eqalign{
&u(\lambda)=u(\lambda_0)-{1\over \pi}\int_{\C}{\pa u(\zeta)\over
\pa\bar\zeta}\biggl({1\over {\zeta-\lambda}}-{1\over {\zeta-\lambda_0}}\biggr)
d\,Re\,\zeta\,d\,Im\,\zeta,\cr
&\lambda\in\C\b 0,\ \ \lambda_0\in\C\b 0,\cr}\eqno(6.4)$$
where $u(\lambda)$ is continuous and bounded for  $\lambda\in\C\b 0$,
$\pa u(\lambda)/\pa\bar\lambda$ is bounded for $\lambda\in\C\b 0$, and
$\pa u(\lambda)/\pa\bar\lambda=O(|\lambda|^{-2})$ as $|\lambda|\to\infty$.
Note that the aforementioned assumptions on $\pa u(\lambda)/\pa\bar\lambda$
in (6.4) can be somewhat weakened. One can prove (6.4) using the formula
$${\pa\over \pa\bar\lambda}{1\over \pi\lambda}=\delta(\lambda) \eqno(6.5)$$
(where $\delta$ is the Dirac function), the Liouville theorem and the
property that (6.4) holds for $\lambda=\lambda_0$.

\vskip 2 mm
{\bf Proposition 6.1.}
{\it Let the assumptions of Lemma} 5.1 {\it be fulfilled}.
{\it Let} $H=H(\lambda,p)$, $R=R(p)$ {\it be defined
by} (6.1), (6.2). {\it Then} $H=H(\lambda,p)$,
$(\lambda,p)\in (\C\b 0)\times (\R^3\b {\cal L}_{\nu})$, {\it satisfies the
following nonlinear integral equation}
$$H(\lambda,p)=R(p)+M(H)(\lambda,p),\ \
\lambda\in\C\b 0,\ \ p\in\R^3\b {\cal L}_{\nu},\eqno(6.6)$$
{\it where}
$$\eqalign{
&M(U)(\lambda,p)=-{1\over \pi}\int_{\C}(U,U)(\zeta,p)
\biggl({1\over {\zeta-\lambda}}-{1\over {\zeta-\lambda_0(p)}}\biggr)
d\,Re\,\zeta\,d\,Im\,\zeta,\cr
&\lambda\in\C\b 0,\ \ p\in\R^3\b {\cal L}_{\nu},\cr}\eqno(6.7)$$
$$\eqalignno{
&(U_1,U_2)(\zeta,p)=\{U_1^{\prime},U_2^{\prime}\}(\zeta,p),\ \
\zeta\in\C\b 0,\ \ p\in\R^3\b {\cal L}_{\nu},&(6.8a)\cr
&U_j^{\prime}(k,p)=U_j(\lambda(k,p),p),\ \ (k,p)\in\Omega_{\nu},\ \ j=1,2,
&(6.8b)\cr}$$
{\it where} $U,U_1,U_2$ {\it are test functions on}
$(\C\b 0)\times (\R^3\b {\cal L}_{\nu})$, $\{U_1^{\prime},U_2^{\prime}\}$
{\it is defined by} (5.5), $\lambda_0=\lambda_0(p)$ {\it is the function of}
(6.2), $\lambda(k,p)$ {\it is defined in} (4.6).

\vskip 2 mm
{\bf Remark 6.1.}
In addition to (6.8), note that definition of $(U_1,U_2)$ can be also
written as
$$\eqalign{
&(U_1,U_2)(\lambda,p)=-{\pi\over 4}\int_{-\pi}^{\pi}\biggl({|p|\over 2}
{{|\lambda|^2-1}\over \bar\lambda |\lambda|}(\cos\v-1)-{|p|\over \bar\lambda}
\sin\v\biggr)\times\cr
&U_1(z_1(\lambda,p,\v),-\xi(\lambda,p,\v))
U_2(z_2(\lambda,p,\v),p+\xi(\lambda,p,\v))d\v,\cr}\eqno(6.9)$$
where
$$\eqalign{
&z_1(\lambda,p,\v)={2k(\lambda,p)(\theta(-\xi(\lambda,p,\v))+i\omega
(-\xi(\lambda,p,\v)))\over i|p|},\cr
&z_2(\lambda,p,\v)={2(k(\lambda,p)+\xi(\lambda,p,\v))
(\theta(p+\xi(\lambda,p,\v))+i\omega
(p+\xi(\lambda,p,\v)))\over i|p|},\cr}\eqno(6.10)$$
$\lambda\in\C\b 0$, $p\in\R^3\b {\cal L}_{\nu}$, $\v\in [-\pi,\pi]$,
$k(\lambda,p)$ is defined in (4.7), $\xi(\lambda,p,\v)$ is defined by (5.2),
$\theta$, $\omega$ are the vector functions of (4.3), (4.4a).

\vskip 2 mm
{\bf Remark 6.2.}
Under the assumptions of Theorem 6.1, equation (6.6) holds, at least, for
almost any $(\lambda,p)\in (\C\b 0)\times (\R^3\b {\cal L}_{\nu})$.

Proposition 6.1 follows from Lemmas 4.1, 5.1, 5.2 and formula (6.4) for
$u(\lambda)=H(\lambda,p)$ (defined by (6.1)).

Consider
$$\eqalign{
&L_{\mu}^{\infty}((\C\b 0)\times (\R^3\b {\cal L}_{\nu}))=\{U\in
L^{\infty}((\C\b 0)\times (\R^3\b {\cal L}_{\nu})):\ |||U|||_{\mu}<
\infty\},\cr
&|||U|||_{\mu}=ess\,\sup\limits_{
\lambda\in\C\b 0,\ p\in\R^3\b {\cal L}_{\nu}}(1+|p|)^{\mu}|U(\lambda,p)|,\ \
\mu>0.\cr}\eqno(6.11)$$
Under the assumptions of Proposition 6.1, from the part I of Proposition 3.1
 and formulas (6.1), (6.2) it follows that
$$H,R\in  L_{\mu}^{\infty}((\C\b 0)\times (\R^3\b {\cal L}_{\nu})) \eqno(6.12)
$$
(where $R$ is independent of $\lambda\in\C\b 0$).

Note that
$$\eqalignno{
&M(U)(\lambda,p)=N(U)(\lambda,p)-N(U)(\lambda_0(p),p),&(6.13a)\cr
&N(U)(\lambda,p)=I(U,U)(\lambda,p),&(6.13b)\cr
&I(U_1,U_2)(\lambda,p)=-{1\over \pi}\int_{\C}(U_1,U_2)(\zeta,p)
{d\,Re\,\zeta\,d\,Im\,\zeta\over {\zeta-\lambda}},&(6.13c)\cr}$$
where  $(\lambda,p)\in (\C\b 0)\times (\R^3\b {\cal L}_{\nu})$, $U$, $U_1$,
$U_2$  are test functions on  $(\C\b 0)\times (\R^3\b {\cal L}_{\nu})$,
$(U_1,U_2)$ is defined by (6.8).

To deal with nonlinear integral equation (6.6) we use Lemmas 6.1, 6.2 and
6.3 given below.

\vskip 2 mm
{\bf Lemma 6.1.}
{\it Let} $U,U_1,U_2\in
L_{\mu}^{\infty}((\C\b 0)\times (\R^3\b {\cal L}_{\nu}))$ {\it for some}
$\mu\ge 2$. {\it Let} $M(U)$, $N(U)$, $I(U_1,U_2)$ {\it be defined by}
(6.7), (6.13), {\it where} $\lambda$, $p$ {\it are the coordinates of
Lemma} 4.1 {\it under assumption} (4.4a). {\it Then}
$$I(U_1,U_2),N(U),M(U)\in
L_{\mu}^{\infty}((\C\b 0)\times (\R^3\b {\cal L}_{\nu})),\eqno(6.14a)$$
$$\eqalign{
&I(U_1,U_2)(\cdot,p),N(U)(\cdot,p),M(U)(\cdot,p)\in
C(\C\b 0)\cap L^{\infty}(\C\b 0)\cr
&{\it for\ almost\ any}\ \ p\in\R^3\b {\cal L}_{\nu},\cr}\eqno(6.14b)$$
$$\eqalignno{
&|||I(U_1,U_2)|||_{\mu}\le c_6(\mu)|||U_1|||_{\mu}|||U_2|||_{\mu},&(6.15a)\cr
&|||N(U)|||_{\mu}\le c_6(\mu)|||U|||_{\mu}^2,&(6.15b)\cr
&|||M(U)|||_{\mu}\le 2c_6(\mu)|||U|||_{\mu}^2,&(6.15c)\cr
&|||N(U_1)-N(U_2)|||_{\mu}\le c_6(\mu)(|||U_1|||_{\mu}+|||U_2|||_{\mu})
|||U_1-U_2|||_{\mu},&(6.16a)\cr
&|||M(U_1)-M(U_2)|||_{\mu}\le 2c_6(\mu)(|||U_1|||_{\mu}+|||U_2|||_{\mu})
|||U_1-U_2|||_{\mu}.&(6.16b)\cr}$$
Lemma 6.1 is proved in Section 11.

\vskip 2 mm
{\bf Lemma 6.2.}
{\it Let} $\mu\ge 2$ {\it and} $0<r<(4c_6(\mu))^{-1}$. {\it Let} $M$ {\it be
defined by} (6.7) ({\it where} $\lambda$, $p$ {\it are the coordinates of
Lemma} 4.1 {\it under assumption} (4.4a)). {\it Let}
$U_0\in L_{\mu}^{\infty}((\C\b 0)\times (\R^3\b {\cal L}_{\nu}))$ {\it and}
$|||U_0|||_{\mu}\le r/2$. {\it Then the equation}
$$U=U_0+M(U) \eqno(6.17)$$
{\it is uniquely solvable for}
$U\in L_{\mu}^{\infty}((\C\b 0)\times (\R^3\b {\cal L}_{\nu}))$,
$|||U|||_{\mu}\le r$, {\it and} $U$ {\it can be found by the method of
successive approximations, in addition},
$$|||U-(M_{U_0})^n(0)|||_{\mu}\le {r(4c_6(\mu)r)^n\over
2(1-4c_6(\mu)r)},\ \ n\in\N,\eqno(6.18)$$
{\it where} $M_{U_0}$ {\it denotes the map} $V\to U_0+M(V)$.

Lemma 6.2 is proved in Section 12 (using Lemma 6.1 and the lemma about
contraction maps).

\vskip 2 mm
{\bf Lemma 6.3.}
{\it Let the assumptions of Lemma} 6.2 {\it be fulfilled}. {\it Let also}
${\tilde U}_0\in L_{\mu}^{\infty}((\C\b 0)\times (\R^3\b {\cal L}_{\nu}))$,
$|||{\tilde U}_0|||_{\mu}\le r/2$ {\it and} $\tilde U$ {\it denote the
solution of} (6.17) {\it with} $U_0$ {\it replaced by} ${\tilde U}_0$,
{\it where}
$\tilde U\in L_{\mu}^{\infty}((\C\b 0)\times (\R^3\b {\cal L}_{\nu}))$,
$|||\tilde U|||_{\mu}\le r$. {\it Then}
$$|||U-\tilde U|||_{\mu}\le {|||U_0-{\tilde U}_0|||_{\mu}\over
{1-4c_6(\mu)r}}.\eqno(6.19)$$
Lemma 6.3 is proved in Section 12.

As a corollary of Proposition 6.1 and Lemmas 6.2 and 6.3, we obtain the
following result.

\vskip 2 mm
{\bf Proposition 6.2.}
{\it Let the assumptions of Lemma} 5.1 {\it be fulfilled. Let}
$$r\buildrel \rm def \over = {2C\over {1-c_1(\mu)C}}
<{1\over 4c_6(\mu)},\eqno(6.20)$$
{\it where} $C$ {\it is the constant of Proposition} 3.1. {\it Let}
$H=H(\lambda,p)$, $R=R(p)$ {\it be defined by}
(6.1), (6.2). {\it Then}
$$|||H|||_{\mu}\le r/2,\ \ |||R|||_{\mu}\le r/2  \eqno(6.21)$$
{\it and} $R$ {\it uniquely and stably determines} $H$ {\it via nonlinear
integral equation} (6.6) {\it considered for} $|||H|||\le r$. {\it In
addition, this equation is solvable by the method of successive approximations
according to } (6.18) ({\it of Lemma} 6.2) {\it and the stability estimate
holds according to}   (6.19) ({\it of Lemma} 6.3) ({\it where} $U_0$,
$U$, ${\tilde U}_0$, $\tilde U$ {\it should be replaced by} $R$, $H$,
$\tilde R$, $\tilde H$, {\it respectively}).

Finally in this section, we apply Propositions 6.1, 6.2 and Lemmas 6.2, 6.3
to approximate finding $H$ on $\Omega$ from $H\big|_{\Xi^{\tau}}$,
where
$$\eqalignno{
&\Omega^{\tau}=\{(k,p)\in\Omega:\ |Im\,k|<\tau\},&(6.22)\cr
&\Xi^{\tau}=\Xi\cap\Omega^{\tau},&(6.23)\cr}$$
where $\Omega$ and $\Xi$ are defined by (2.19), (2.20).
In the coordinates of  Lemma 4.1
this means that we deals with approximate finding $H=H(\lambda,p)$
defined by (6.1) from $R_{2\tau}=\chi_{2\tau}R$, where $R=R(p)$ is defined by
 (6.2) and
$\chi_s$ denotes the multiplication operator by the function
$\chi_r(p)$, where
$$\chi_s(p)=1\ \ {\rm for}\ \ |p|<s,\ \ \chi_s(p)=0\ \ {\rm for}\ \ |p|\ge s,
\ \ {\rm where}\ \ p\in\R^3,\ s>0.\eqno(6.24)$$
One can see that $R_{2\tau}$ is a low-frequency part of $R$ and, thus,
$H\big|_{\Xi^{\tau}}$ is a low-frequency part of
$H\big|_{\Xi}$. One can see also that $\Omega^{\tau}$
is a low-imaginary part of $\Omega$ and, therefore, $\Xi^{\tau}$
is a low-imaginary part of $\Xi$.

Note that
$$\eqalignno{
&|||\chi_{2\tau}R|||_{\mu}\le |||R|||_{\mu},&(6.25)\cr
&|||R-\chi_{2\tau}R|||_{\mu}\le {|||R|||_{\mu^*}\over (1+2\tau)^{\mu^*-\mu}}
&(6.26)\cr}$$
for $R\in L_{\mu^*}^{\infty}((\C\b 0)\times (\R^3\b {\cal L}_{\nu}))$, where
$0\le\mu\le\mu^*$, $\tau>0$.

Using Propositions 6.1, 6.2, Lemmas 6.2, 6.3 and estimates (6.25), (6.26)
we obtain the following result.

\vskip 2 mm
{\bf Proposition 6.3.}
{\it Let the assumptions of Proposition} 6.2 {\it be fulfilled.
Let also}
$$\hat v\in L_{\mu^*}^{\infty}(\R^3)\ \ {\it for\ some}\ \ \mu^*>\mu.
\eqno(6.27)$$
{\it Let} $\tau>0$. {\it Then}:
$$\eqalignno{
&|||\chi_{2\tau}R|||_{\mu}\le r/2,&(6.28a)\cr
&R\in L_{\mu^*}^{\infty}(\R^3\b {\cal L}_{\nu});&(6.28b)\cr}$$
$\chi_{2\tau}R$  {\it uniquely and stably determines} $H_{2\tau}$,
{\it where}  $H_{2\tau}$ {\it denotes the solution of  the
nonlinear integral equation}
$$H_{2\tau}=\chi_{2\tau}R+M(H_{2\tau}),\ \  |||H_{2\tau}|||_{\mu}\le r,
\eqno(6.29)$$
{\it see Lemmas} 6.2, 6.3; {\it the following estimate holds}:
$$|||H-H_{2\tau}|||_{\mu}\le {|||R|||_{\mu^*}\over (1+2\tau)^{\mu^*-\mu}
(1-4c_6(\mu)r)}.\eqno(6.30)$$
Note that (6.28b) follows from the property that
$R\in L_{\mu}^{\infty}(\R^3\b {\cal L}_{\nu})$, the assumption (6.27), the
part II of Proposition 3.1 for $\mu=\mu^*$ and definition (6.2).
Estimate (6.30) follows from Proposition 6.2, Lemma 6.3 (where $U_0$, $U$,
${\tilde U}_0$, $\tilde U$ are replaced by $R$, $H$, $\chi_{2r}R$,
$H_{2\tau}$, respectively) and from (6.28), (6.29), (6.26).

Actually, in Proposition 6.3, $H_{2\tau}$ is a low-frequency approximation
to $H$. In addition, estimate (6.30) shows that the error between
$H_{2\tau}$ and $H$ rapidly decays in the norm $|||\cdot|||_{\mu}$ as
$\tau\to +\infty$ if $\mu^*-\mu$ is sufficiently great.

\vskip 2 mm
{\bf 7. Finding $\hat v$ on $\R^3$ from $H$ on $\Omega$ and some related
results}
\vskip 2 mm
Actually, in this section we consider finding
$\hat v$ on $\R^3\b {\cal L}_{\nu}$ from $H$ on $\Omega_{\nu}$
in the coordinates of
Lemma 4.1 under assumption (4.4a). In addition, under the assumptions of
Proposition 6.3, we consider also approximate finding
$\hat v$ on $\R^3\b {\cal L}_{\nu}$ from $H_{2\tau}$ introduced in
Proposition 6.3 as a low-frequency approximation to $H$.

Under assumption (2.3), formulas (2.11), (4.7), (4.8) imply that
$$\eqalignno{
&H(\lambda,p)\to \hat v(p)\ \ {\rm as}\ \ \lambda\to 0,&(7.1a)\cr
&H(\lambda,p)\to \hat v(p)\ \ {\rm as}\ \ \lambda\to \infty,&(7.1b)\cr}$$
where $\lambda\in\C\b 0$, $p\in\R^3\b {\cal L}_{\nu}$ and $H(\lambda,p)$ is
defined by (6.1). In addition, under the assumptions of Proposition 6.1,
formulas (6.6), (7.1) (and estimates (3.7), (3.8), (5.7), (5.8)) imply that
$$\eqalignno{
&\hat v(p)=R(p)+M(H)(0,p),&(7.2a)\cr
&\hat v(p)=R(p)-N(H)(\lambda_0(p),p)&(7.2b)\cr}$$
for $p\in\R^3\b {\cal L}_{\nu}$, where $M$, $N$ are defined by (6.7), (6.8),
(6.13), and $\lambda_0$ is the function of (6.2). In addition, due to
(6.13a), we have that
$$M(H)(0,p)=N(H)(0,p)-N(\lambda_0(p),p),\ \ p\in\R^3\b {\cal L}_{\nu},
\eqno(7.3)$$
and, as a corollary of (7.2), (7.3), we have that
$$N(H)(0,p)=0,\ \  p\in\R^3\b {\cal L}_{\nu}.\eqno(7.4)$$

Further, under the assumptions of Proposition 6.3, using (6.29) we obtain
that
$$\eqalignno{
&H_{2\tau}(\lambda,p)\to \hat v^+_{2\tau}(p)\ \ {\rm as}\ \ \lambda\to 0,
&(7.5a)\cr
&H_{2\tau}(\lambda,p)\to \hat v^-_{2\tau}(p)\ \ {\rm as}\ \ \lambda\to\infty,
&(7.5b)\cr}$$
where
$$\eqalignno{
&\hat v^+_{2\tau}(p)=\chi_{2\tau} R(p)+M(H_{2\tau})(0,p),&(7.6a)\cr
&\hat v^-_{2\tau}(p)=\chi_{2\tau} R(p)-N(H_{2\tau})(\lambda_0(p),p),&(7.6b)
\cr}$$
for $p\in\R^3\b {\cal L}_{\nu}$, where $M$, $N$ are defined by (6.7), (6.8),
(6.13) and $\lambda_0$ is the function of (6.2). In addition, formulas
(1.9), (6.11), (7.1), (7.5) imply that
$$\|\hat v-\hat v^{\pm}_{2\tau}\|_{\mu}\le |||H-H_{2\tau}|||_{\mu}.
\eqno(7.7)$$
Under the assumptions of Proposition 6.3, formulas (6.30), (7.7) imply that
$\hat v$ on $\R^3$ can be approximately determined from $H_{2\tau}$ as
$\hat v^{\pm}_{2\tau}$ of (7.5), (7.6) and
$$\|\hat v-\hat v^{\pm}_{2\tau}\|_{\mu}=O\bigl({1\over \tau^{\mu^*-\mu}}
\bigr)\ \ {\rm as}\ \ \tau\to +\infty.\eqno(7.8)$$

\vskip 2 mm
{\bf 8. Proofs of Lemma 3.1 and Proposition 3.1}

{\it Proof of} (3.3). We have that
$$|A(k)U(p)|\le I(k,p)\|\hat v\|_{\mu}\|U\|_{\mu},\eqno(8.1)$$
where
$$I(k,p)=\int\limits_{\R^3}
{d\xi\over (1+|p+\xi|)^{\mu}(1+|\xi|)^{\mu}
|\xi^2+2k\xi|},\ \ k\in\Sigma,\ \ p\in\R^3.\eqno(8.2)$$
To prove (3.3) it is sufficient to prove that
$$\eqalignno{
&I(k,p)\le {c_1(\mu)\over (1+|p|)^{\mu}},&(8.3a)\cr
&I(k,p)\le {c_2(\mu)(\ln\,(|Im\,k|))^2\over |Im\,k|(1+|p|)^{\mu}},\ \
\ln\,|Im\,k|\ge 2,&(8.3b)\cr}$$
where $k\in\Sigma$, $p\in\R^3$. Note that
$$I(k,p)\le \biggl(\int\limits_{|\xi|\le |p+\xi|}+
\int\limits_{|\xi|\ge |p+\xi|}\biggr)
{d\xi\over (1+|p+\xi|)^{\mu}(1+|\xi|)^{\mu}|\xi^2+2k\xi|},\eqno(8.4)$$
where $k\in\Sigma$, $p\in\R^3$. Note also that
$$|\xi|\le |p+\xi|\Rightarrow |p+\xi|\ge |p|/2,\ \
|\xi|\ge |p+\xi|\Rightarrow |\xi|\ge |p|/2,\eqno(8.5)$$
where $\xi,p\in\R^3$. Using (8.4), (8.5) we obtain that
$$I(k,p)\le (1+|p|/2)^{-\mu}(I_1(k)+I_2(k,p)),\eqno(8.6)$$
where
$$\eqalign{
&I_1(k)=\int\limits_{\R^3}{d\xi\over (1+|\xi|)^{\mu}|\xi^2+2k\xi|},\cr
&I_2(k,p)=\int\limits_{\R^3}{d\xi\over (1+|p+\xi|)^{\mu}|\xi^2+2k\xi|},
\cr}\eqno(8.7)$$
where $k\in\Sigma$, $p\in\R^3$. Note that
$$I_1(k)=I_2(k,0),\ \ k\in\Sigma.\eqno(8.8)$$
Note further that
$$\eqalign{
&I_2(k,p)=\cr
&\int\limits_{\R^3}{d\xi\over (1+|(\xi+Re\,k)-(Re\,k-p)|)^{\mu}
|(\xi+Re\,k)^2-(Re\,k)^2+2iIm\,k(\xi+Re\,k)|}=\cr
&I_3(k,Re\,k-p),\cr}\eqno(8.9)$$
$$I_3(k,p)=
\int\limits_{\R^3}{d\xi\over (1+|\xi-p|)^{\mu}
|\xi^2-(Re\,k)^2+2iIm\,k\xi|},\eqno(8.10)$$
where $k\in\Sigma$, $p\in\R^3$.

In view of (8.6)-(8.10), to prove (8.3) it is sufficient to prove that
$$\eqalignno{
&I_3(k,p)\le \tilde c_1(\mu),&(8.11a)\cr
&I_3(k,p)\le {\tilde c_2(\mu)(\ln\,(|Im\,k|))^2\over |Im\,k|},\ \
|Im\,k|\ge \rho(\mu),&(8.11b)\cr}$$
where  $k\in\Sigma$, $p\in\R^3$.

Consider $p_{\pr}=p_{\pr}(p,Im\,k)$, $p_{\perp}=p_{\perp}(p,Im\,k)$,
where
$$p_{\pr}={p\,Im\,k\over |Im\,k|} {Im\,k\over |Im\,k|}\ \ {\rm for}\ \
|Im\,k|\ne 0,\ \ p_{\pr}=0\ \ {\rm for}\ \ |Im\,k|=0,\ p_{\perp}=p-p_{\pr},
\eqno(8.12)$$
where $p,Im\,k\in\R^3$. Using the properties
$$Im\,k\,Re\,k=0,\ \ (Re\,k)^2=(Im\,k)^2\ \ {\rm for}\ \ k\in\Sigma
\eqno(8.13)$$
and changing variables in the integral of (8.10), we obtain that
$$\eqalign{
&I_3(k,p)=\cr
&\int\limits_{\R^3}{d\xi\over
(1+((\xi_1-|p_{\perp}|)^2+\xi_2^2+(\xi_3-sgn\,(p_{\pr}Im\,k)|p_{\pr}|)^2
)^{1/2})^{\mu}|\xi^2-(Im\,k)^2+2i\,|Im\,k|\xi_3|},\cr}\eqno(8.14)$$
where $k\in\Sigma$, $p\in\R^3$. Further, using (8.14) we obtain that
$$I_3(k,p)\le \sqrt{2}I_4(|Im\,k|,|p_{\pr}|,|p_{\perp}|),\eqno(8.15)$$
$$\eqalign{
&I_4(\rho,s,t)=\cr
&\int\limits_{\R^3}{d\xi\over (1+(\xi_1-t)^2+\xi_2^2+(\xi_3-s)^2)^{\mu/2}
(|\xi_1^2+\xi_2^2+\xi_3^2-\rho^2|+2\rho |\xi_3|)},\cr}\eqno(8.16)$$
where $k\in\Sigma$, $p\in\R^3$, $\rho,s,t\in [0,+\infty [$.
Due to (8.15), to prove (8.11) it is sufficient to prove that

$$\eqalignno{
&I_4(\rho,s,t)\le \tilde c_1(\mu)/\sqrt{2},&(8.17a)\cr
&I_4(\rho,s,t)\le {\tilde c_2(\mu)(\ln\rho)^2\over \sqrt{2}\rho},\ \
\ln\rho\ge 2,&(8.17b)\cr}$$
where $\rho,s,t\in [0,+\infty [$.  Note that
$$\eqalign{
&I_4(\rho,s,t)\le \biggl(\int\limits_{|\xi_3|\le |\xi_3-s|}+
\int\limits_{|\xi_3|\ge |\xi_3-s|}\biggr)\times\cr
&{d\xi\over (1+(\xi_1-t)^2+\xi_2^2+(\xi_3-s)^2)^{\mu/2}
(|\xi_1^2+\xi_2^2+\xi_3^2-\rho^2|+2\rho |\xi_3|)}\le\cr
&\int\limits_{|\xi_3|\le |\xi_3-s|}
{d\xi\over (1+(\xi_1-t)^2+\xi_2^2+\xi_3^2)^{\mu/2}
(|\xi_1^2+\xi_2^2+\xi_3^2-\rho^2|+2\rho |\xi_3|)}+\cr
&\int\limits_{|\xi_3|\ge |\xi_3-s|}
{d\xi\over (1+(\xi_1-t)^2+\xi_2^2+(\xi_3-s)^2)^{\mu/2}
(|\xi_1^2+\xi_2^2+(\xi_3-s)^2-\rho^2|+2\rho |\xi_3-s|)}\le\cr
&2\int\limits_{\R^3}
{d\xi\over (1+(\xi_1-t)^2+\xi_2^2+\xi_3^2)^{\mu/2}
(|\xi_1^2+\xi_2^2+\xi_3^2-\rho^2|+2\rho |\xi_3|)}=2I_4(\rho,0,t),\cr}
\eqno(8.18)$$
where $\rho,s,t\in [0,+\infty [$. In addition, in (8.18) we used, in
particular, that
$$|\xi_1^2+\xi_2^2+\xi_3^2-\rho^2|+2\rho |\xi_3|\ge
|\xi_1^2+\xi_2^2+(\xi_3-s)^2-\rho^2|+2\rho |\xi_3-s|\ \ {\rm if}\ \
|\xi_3|\ge |\xi_3-s|.\eqno(8.19)$$
To prove (8.19) we rewrite it as
$$\eqalign{
&\rho^2-\xi_1^2-\xi_2^2-\xi_3^2+2\rho |\xi_3|\ge
\rho^2-\xi_1^2-\xi_2^2-(\xi_3-s)^2+2\rho |\xi_3-s|\cr
&{\rm for}\ \ |\xi_3|\ge |\xi_3-s|,\ \ \xi_1^2+\xi_2^2+\xi_3^2\le\rho^2,\cr}
\eqno(8.20a)$$
$$\eqalign{
&\xi_1^2+\xi_2^2+\xi_3^2-\rho^2+2\rho |\xi_3|\ge
\xi_1^2+\xi_2^2+(\xi_3-s)^2-\rho^2+2\rho |\xi_3-s|\cr
&{\rm for}\ \ |\xi_3|\ge |\xi_3-s|,\ \ \xi_1^2+\xi_2^2+\xi_3^2\ge\rho^2,\ \
\xi_1^2+\xi_2^2+(\xi_3-s)^2\ge\rho^2,\cr}\eqno(8.20b)$$
$$\eqalign{
&\xi_1^2+\xi_2^2+\xi_3^2-\rho^2+2\rho |\xi_3|\ge
\rho^2-\xi_1^2-\xi_2^2-(\xi_3-s)^2+2\rho |\xi_3-s|\cr
&{\rm for}\ \ |\xi_3|\ge |\xi_3-s|,\ \ \xi_1^2+\xi_2^2+\xi_3^2\ge\rho^2,\ \
\xi_1^2+\xi_2^2+(\xi_3-s)^2\le\rho^2.\cr}\eqno(8.20c)$$
Inequality (8.20a) follows from the inequalities
$$\eqalignno{
&-x^2+2\rho x\ge -y^2+2\rho y\ \ {\rm for}\ \ 0\le y\le x\le\rho,&(8.21)\cr
&y=|\xi_3-s|\le x=|\xi_3|\le\sqrt{\rho^2-\xi_1^2+\xi_2^2}\le\rho.&(8.22)\cr}$$
Inequality (8.20b) is obvious. Inequality (8.20c) follows from the
inequalities
$$\eqalignno{
&x^2-\delta^2+2\rho x\ge\delta^2-y^2+2\rho y\ \ {\rm for}\ \
0\le\delta\le\rho,\ \ 0\le y\le \delta\le x,&(8.23)\cr
&y=|\xi_3-s|\le \delta=\sqrt{\rho^2-\xi_1^2-\xi_2^2}\le x=|\xi_3|,\ \
\delta=\sqrt{\rho^2-\xi_1^2-\xi_2^2}\le\rho.&(8.24)\cr}$$
In turn, inequality (8.23) follows from the inequalities
$$\eqalign{
&x^2-\delta^2+2\rho x\ge 2\rho\delta\ \ {\rm for}\ \ 0\le\delta\le x,\cr
&\delta^2-y^2+2\rho y\buildrel (8.21) \over \le 2\rho\delta\ \ {\rm for}\ \
0\le y\le\delta\le\rho.\cr}\eqno(8.25)$$
Thus formulas (8.19), (8.18) are proved.

Due to (8.18), to prove (8.17) it is sufficient to prove that
$$\eqalignno{
&I_4(\rho,0,t)\le {\tilde c_1(\mu)\over 2\sqrt{2}},&(8.26a)\cr
&I_4(\rho,0,t)\le {\tilde c_2(\mu)(\ln\rho)^2\over 2\sqrt{2}\rho},\ \
\ln\rho\ge 2,&(8.26b)\cr}$$
where $\rho,t\in [0,+\infty [$. Using spherical coordinates we obtain that
$$\eqalign{
&I_4(\rho,0,t)=\cr
&\int\limits_0^{+\infty}\int\limits_{-\pi}^{\pi}\int\limits_0^{\pi}
{r^2\sin\psi d\psi d\v dr\over
(1+r^2+t^2-2rt\sin\psi \cos\v)^{\mu/2}(|r^2-\rho^2|+2\rho\,r|\cos\psi|)}=\cr
&2\int\limits_0^{+\infty}\int\limits_{-\pi}^{\pi}\int\limits_0^{\pi/2}
{r^2\sin\psi d\psi d\v dr\over
(1+r^2+t^2-2rt\sin\psi \cos\v)^{\mu/2}(|r^2-\rho^2|+2\rho\,r\cos\psi)}\le\cr
&4\int\limits_0^{+\infty}\int\limits_{-\pi/2}^{\pi/2}\int\limits_0^{\pi/2}
{r^2\sin\psi d\psi d\v dr\over
(1+r^2+t^2-2rt\sin\psi \cos\v)^{\mu/2}(|r^2-\rho^2|+2\rho\,r\cos\psi)}\le\cr
&4\int\limits_0^{+\infty}\int\limits_{-\pi/2}^{\pi/2}\int\limits_0^{\pi/2}
{r^2\sin\psi d\psi d\v dr\over
(1+r^2+t^2-2rt\cos\v)^{\mu/2}(|r^2-\rho^2|+2\rho\,r\cos\psi)}=\cr
&4\int\limits_0^{+\infty}\biggl(\int\limits_{-\pi/2}^{\pi/2}
{d\v\over(1+r^2+t^2-2rt\cos\v)^{\mu/2}}
\int\limits_0^{\pi/2}
{\sin\psi d\psi\over(|r^2-\rho^2|+2\rho\,r\cos\psi)}\biggr)r^2dr,\cr}
\eqno(8.27)$$
where $\rho,t\in [0,+\infty [$. Further, we obtain that:
$$\eqalign{
&\int\limits_{-\pi/2}^{\pi/2}{d\v\over {1+r^2+t^2-2rt\cos\v}}=
\int\limits_{-\pi/2}^{\pi/2}{d\v\over {1+r^2+t^2-2rt(1-2(\sin(\v/2))^2)}}=\cr
&2\int\limits_{-\pi/4}^{\pi/4}{d\v\over {1+(r-t)^2+4rt(\sin\v)^2}}\le
2\int\limits_{-\pi/4}^{\pi/4}
{\sqrt{2}\cos\v d\v\over {1+(r-t)^2+4rt(\sin\v)^2}}=\cr
&4\int\limits_0^{1/\sqrt{2}}{\sqrt{2}du\over {1+(r-t)^2+4rtu^2}}=
{2\sqrt{2}\over \sqrt{rt}}\int\limits_0^{\sqrt{2rt}}{du\over {1+(r-t)^2+u^2}}
\le\cr
&{4\sqrt{2}\over \sqrt{rt}}\int\limits_0^{\sqrt{2rt}}
{du\over (\sqrt{1+(r-t)^2}+u)^2}=\cr
&{4\sqrt{2}\over \sqrt{rt}}\biggl({1\over \sqrt{1+(r-t)^2}}-
{1\over {\sqrt{1+(r-t)^2}+\sqrt{2rt}}}\biggr)=\cr
&{8\over \sqrt{1+(r-t)^2}(\sqrt{1+(r-t)^2}+\sqrt{2rt})};\cr}\eqno(8.28)$$
$$\eqalign{
&\int\limits_0^{\pi/2}{\sin\psi d\psi\over {|r^2-\rho^2|+2\rho\,r\cos\psi}}=
\int\limits_0^1{du\over {|r^2-\rho^2|+2\rho\,ru}}=\cr
&{1\over 2\rho\,r}\ln\,\bigl(|r^2-\rho^2|+2\rho\,ru\bigr)\big|_0^1=
{1\over 2\rho\,r}\ln\,\biggl(1+{2\rho\,r\over |r^2-\rho^2|}\biggr),\cr}
\eqno(8.29)$$
where $\rho,t\in [0,+\infty [$. Using (8.27)-(8.29) we obtain that
$$\eqalign{
&I_4(\rho,0,t)\buildrel \mu\ge 2 \over \le 32\int\limits_0^{+\infty}
{dr\over  \sqrt{1+(r-t)^2}(\sqrt{1+(r-t)^2}+\sqrt{2rt})}\le\cr
&32\int\limits_0^{+\infty}{dr\over {1+(r-t)^2}}\le
32\int\limits_{-\infty}^{+\infty}{dr\over {1+r^2}}=32\pi\ \ {\rm for}\ \
\rho=0,\ \ t\ge 0,\cr}\eqno(8.30a)$$

$$\eqalign{
&I_4(\rho,0,t)\buildrel \mu\ge 2 \over \le {16\over \rho}
\int\limits_0^{+\infty}{\ln\,\biggl(1+{2\rho r\over |r^2-\rho^2|}\biggr)rdr
\over  \sqrt{1+(r-t)^2}(\sqrt{1+(r-t)^2}+\sqrt{2rt})}=\cr
&I_5(\rho,t/\rho)\ \ {\rm for}\ \ \rho>0,\ \ t\ge 0,\cr}\eqno(8.30b)$$
where
$$
I_5(\rho,\ep)=16
\int\limits_0^{+\infty}{\rho\ln\,\biggl(1+{2\tau\over |\tau^2-1|}\biggr)
\tau d\tau
\over  \sqrt{1+\rho^2(\tau-\ep)^2}(\sqrt{1+\rho^2(\tau-\ep)^2}+
\rho\sqrt{2\tau\ep})},\ \ \rho>0,\ \ep\ge 0.$$

As regards $I_5(\rho,\ep)$, we will estimate it separately for
$\ep\in [0,1/4]$, $\ep\in [1/4,2]$ and $\ep\in [2,+\infty [$. For
$\ep\in [0,1/4]$, $\rho>0$, we start with the partition:
$$\eqalign{
&I_5(\rho,\ep)=16\biggl(\int\limits_0^{1/2}+\int\limits_{1/2}^{3/2}+
\int\limits_{3/2}^{+\infty}\biggr)
{\rho\ln\,\biggl(1+{2r\over |r^2-1|}\biggr)rdr\over
{1+\rho^2(r-\ep)^2+\rho\sqrt{2r\ep}\sqrt{1+\rho^2(r-\ep)^2}}}=\cr
&16(I_{5,1}(\rho,\ep)+I_{5,2}(\rho,\ep)+I_{5,3}(\rho,\ep)),\cr}\eqno(8.31)$$
where $I_{5,1},I_{5,2},I_{5,3}$ correspond to
$\int\limits_0^{1/2}$, $\int\limits_{1/2}^{3/2}$,
$\int\limits_{3/2}^{+\infty}$, respectively. Further,
$$I_{5,1}(\rho,\ep)\le\ln\,(7/3)\int\limits_0^{1/2}{\rho\,rdr\over
{1+\rho^2(r-\ep)^2+\rho\sqrt{r\ep}(1+\rho |r-\ep|)}}=\ln\,(7/3)\tilde I_{5,1}
(\rho,\ep),\eqno(8.32)$$
where $\rho>0$, $\ep\in [0,1/4]$. In addition:
$$\tilde I_{5,1}(\rho,\ep)=\int\limits_0^{1/2}{\rho\,rdr\over {1+\rho^2r^2}}=
{1\over 2}\int\limits_0^{1/4}{\rho\,d\tau\over {1+\rho^2\tau}}=
{\ln\,\bigl(1+\rho^2/4\bigr)\over 2\rho}\eqno(8.33a)$$
for $\rho>0$, $\ep=0$;
$$\eqalign{
&\tilde I_{5,1}(\rho,\ep)=\biggl(\int\limits_0^{1/2}+\int\limits_{1/2}^{3/2}+
\int\limits_{3/2}^{1/(2\ep)}\biggr)
{\rho\,\ep^2\tau d\tau\over
{1+(\rho\ep)^2(\tau-1)^2+\rho\ep\sqrt{\tau}(1+\rho\ep|\tau-1|)}}=\cr
&\tilde I_{5,1,1}(\rho,\ep)+\tilde I_{5,1,2}(\rho,\ep)+\tilde I_{5,1,\ep}
(\rho,\ep)\ \ {\rm for}\ \ \rho>0,\ \ \ep\in ]0,1/4],\cr}\eqno(8.33b)$$
where
$\tilde I_{5,1,1}$, $\tilde I_{5,1,2}$, $\tilde I_{5,1,\ep}$ correspond
to
$\int\limits_0^{1/2}$, $\int\limits_{1/2}^{3/2}$,
$\int\limits_{3/2}^{1/(2\ep)}$, respectively. In addition:
$$\tilde I_{5,1,1}(\rho,\ep)\le {\rho\ep^2/4\over {1+(\rho\ep)^2/4}}=
{(\rho\ep)^2\over \rho\,(4+(\rho\ep)^2)}\le \min\,\bigl({1\over \rho},
{\rho\over 4^3}\bigr),\eqno(8.34)$$
$$\eqalign{
&\tilde I_{5,1,2}(\rho,\ep)\le 2\int\limits_1^{3/2}
{\rho\,\ep^2(3/2) d\tau\over
{1+\rho\ep\sqrt{1/2}(1+\rho\ep(\tau-1))}}=\cr
&3\int\limits_0^{1/2}
{\rho\,\ep^2 d\tau\over
{1+\rho\ep/\sqrt{2}+(\rho\ep)^2\tau/\sqrt{2}}}=
{3\sqrt{2}\rho\ep^2\over (\rho\ep)^2}\ln\,(\sqrt{2}+\rho\ep+(\rho\ep)^2\tau)
\big|_0^{1/2}=\cr
&{3\sqrt{2}\over \rho}\ln\,\bigl(1+{(\rho\ep)^2\over 2(\sqrt{2}+\rho\ep)}
\bigr)\le
{3\sqrt{2}\over \rho}\ln\,\bigl(1+{\rho^2\over 32\sqrt{2}}\bigr),\cr}
\eqno(8.35)$$
$$\eqalign{
&\tilde I_{5,1,3}(\rho,\ep)\le \int\limits_{1/2}^{(2\ep)^{-1}-1}
{\rho\,\ep^2 (\tau+1) d\tau\over{1+(\rho\ep)^2\tau^2}}\le
3\int\limits_{1/2}^{(2\ep)^{-1}}
{\rho\,\ep^2 \tau d\tau\over{1+(\rho\ep)^2\tau^2}}=\cr
&{3\rho\ep^2\over 2(\rho\ep)^2}\ln\,(1+(\rho\ep)^2x)\big|_{1/4}^{1/(2\ep)^2}\le
{3\over 2\rho}\ln\,\bigl(1+\rho^2/4\bigr),\cr}\eqno(8.36)$$
where $\rho>0$, $\ep\in ]0, 1/4]$. Further,
$$\eqalignno{
&I_{5,2}(\rho,\ep)\le {\rho\over {1+\rho^2/16}}\int\limits_{1/2}^{3/2}
\ln\,\bigl(1+{2r\over |r^2-1|}\bigr)rdr,&(8.37)\cr
&I_{5,3}(\rho,\ep)\le \int\limits_{3/2}^{+\infty}
{\rho(2r/|r^2-1|)rdr\over {1+\rho^2(r-\ep)^2}}\le
\int\limits_{3/2}^{+\infty}
{4\rho(1+(r^2-1)^{-1})dr\over (1+\rho\,(r-\ep))^2}\le {8\over {1+\rho\,5/4}},
&(8.38)\cr}$$
where  $\rho>0$, $\ep\in ]0, 1/4]$.

For $\ep\in [1/4,2]$, $\rho>0$ we use the partition:
$$\eqalign{
&I_5(\rho,\ep)=16\biggl(\int\limits_0^{1/8}+\int\limits_{1/8}^3+
\int\limits_3^{+\infty}\biggr)
{\rho\ln\,\biggl(1+{2r\over |r^2-1|}\biggr)rdr\over
{1+\rho^2(r-\ep)^2+\rho\sqrt{2r\ep}\sqrt{1+\rho^2(r-\ep)^2}}}=\cr
&16(I_{5,4}(\rho,\ep)+I_{5,5}(\rho,\ep)+I_{5,6}(\rho,\ep)),\cr}\eqno(8.39)$$
where $I_{5,4},I_{5,5},I_{5,6}$ correspond to
$\int\limits_0^{1/8}$, $\int\limits_{1/8}^3$,
$\int\limits_3^{+\infty}$, respectively. In addition:
$$\eqalignno{
&I_{5,4}(\rho,\ep)\le\ln\,(3/2)\int\limits_0^{1/8}{\rho\,rdr\over
{1+\rho^2(r-\ep)^2}}\le {\ln\,(3/2)\rho\over {64+\rho^2}},&(8.40)\cr
&I_{5,5}(\rho,\ep)\le \rho\int\limits_{1/8}^3
\ln\,\bigl(1+{2r\over |r^2-1|}\bigr)rdr,&(8.41a)\cr}$$
$$\eqalign{
&I_{5,5}(\rho,\ep)\le\int\limits_{1/8}^3
{3\rho\ln\,\biggl(1+{2\over |r-1|}\biggr)dr\over
{1+\rho\sqrt{1/32}(1+\rho\,|r-\ep|)}}\le\cr
&\biggl(\int\limits_{|r-\ep|\le |r-1|,\,1/8<r<3}+
\int\limits_{|r-\ep|\ge |r-1|,\,1/8<r<3}\biggr)
{3\rho\ln\,\biggl(1+{2\over |r-1|}\biggr)dr\over
{1+\rho\sqrt{1/32}+\rho^2\sqrt{1/32}|r-\ep|)}}\le\cr
&12\sqrt{2}\biggl(\int\limits_{1/8}^3
{\rho\ln\,\biggl(1+{2\over |r-\ep|}\biggr)dr\over
{4\sqrt{2}+\rho+\rho^2|r-\ep|}}+
\int\limits_{1/8}^3
{\rho\ln\,\biggl(1+{2\over |r-1|}\biggr)dr\over
{4\sqrt{2}+\rho+\rho^2|r-1|}}\biggr)\le\cr
&48\sqrt{2}\int\limits_0^3
{\rho\ln\,\bigl(1+{2\over r}\bigr)dr\over
{4\sqrt{2}+\rho+\rho^2r}}\le
48\sqrt{2}\int\limits_0^3
{\ln\,\bigl(1+{2\over r}\bigr)dr\over
{1+\rho\,r}}=\cr
&{48\sqrt{2}\over \rho}\biggl(\int\limits_0^1+\int\limits_1^{3\rho}\biggr)
{\ln\,\bigl(1+{2\rho\over \tau}\bigr)d\tau\over
{1+\tau}}\buildrel \rho\ge 1 \over \le\cr
&{48\sqrt{2}\over \rho}\biggl(\int\limits_0^1
{(\ln\,(3\rho)+\ln\,(1/\tau))d\tau\over {1+\tau}}
+\int\limits_1^{3\rho}
{\ln\,(1+2\rho)d\tau\over{1+\tau}}\biggr)=\cr
&{48\sqrt{2}\over \rho}\biggl(\ln\,(1+2\rho)\,\ln\,(1/2+(3/2)\rho)+
\ln\,(3\rho)\,\ln\,2+\int\limits_0^1
{\ln\,(1/\tau)\,d\tau\over {1+\tau}}\biggr),\ \ \rho\ge 1,\cr}\eqno(8.41b)$$
$$
I_{5,6}(\rho,\ep)\le \int\limits_3^{+\infty}
{\rho(2r/|r^2-1|)rdr\over {1+\rho^2(r-\ep)^2}}\le
\int\limits_3^{+\infty}
{4\rho(1+(r^2-1)^{-1})dr\over (1+\rho\,(r-\ep))^2}\le {5\over {1+\rho}},
\eqno(8.42)$$
where $\rho>0$, $\ep\in [1/4,2]$.

For $\ep\in [2,+\infty [$, $\rho>0$ we use the partition:
$$I_5(\rho,\ep)=16(I_{5,7}(\rho,\ep)+I_{5,8}(\rho,\ep)),\eqno(8.43)$$
where $I_{5,7}=I_{5,1}+I_{5,2}$,
 $I_{5,8}=I_{5,3}$, where $I_{5,1}$,
$I_{5,2}$, $I_{5,3}$ are defined as in (8.31). In addition,
$$I_{5,7}(\rho,\ep)\le {\rho\over {1+\rho^2/4}}\int\limits_0^{3/2}
\ln\,\bigl(1+{2r\over |r^2-1|}\bigr)rdr,\eqno(8.44)$$
$$\eqalign{
&I_{5,8}(\rho,\ep)\le\int\limits_{3/2}^{+\infty}
{\rho(2r/|r^2-1|)rdr\over
{1+\rho^2(r-\ep)^2+\rho\sqrt{3}(1+\rho\,|r-\ep|)}}=\cr
&\int\limits_{3/2}^{+\infty}
{2\rho(1+(r^2-1)^{-1})dr\over
{1+\sqrt{3}\rho+\sqrt{3}\rho^2|r-\ep|+\rho^2(r-\ep)^2}}\le\cr
&\int\limits_{-\infty}^{+\infty}
{4\rho\,dr\over
{1+\sqrt{3}\rho+\sqrt{3}\rho^2|r|+\rho^2r^2}}\le
\int\limits_0^1{8\rho\,dr\over {\sqrt{3}\rho+\sqrt{3}\rho^2r}}+
\int\limits_1^{+\infty}{8\rho\,dr\over {1+\rho^2r^2}}\le\cr
&{8\over \sqrt{3}}\int\limits_0^1{dr\over {1+\rho\,r}}+
\int\limits_1^{+\infty}{16\rho\,dr\over (1+\rho\,r)^2}=
{8\over \sqrt{3}} {\ln\,(1+\rho)\over \rho}+{16\over {1+\rho}},\cr}\eqno(8.45)
$$
where $\rho>0$, $\ep\in [2,+\infty [$.

Estimates (8.26) follow from (8.30)-(8.45). Thus, estimates (8.17), (8.11),
(8.3) are proved. The proof of (3.3) is completed.

\vskip 2 mm
{\it Proof of} (3.2).
Let
$$f_1(\xi)=\hat v(\xi),\ \ f_2(\xi)={U(\xi)\over {\xi^2-2k\xi}},\eqno(8.46)$$
where $\xi\in\R^3$, $k\in\Sigma$. We have, in particular, that
$$f_1\in L^{\infty}(\R^3),\ \  f_2\in L^1(\R^3).\eqno(8.47)$$
Property (3.2) follows from (8.46), (8.47) and the following lemma.

\vskip 2 mm
{\bf Lemma 8.1.}
{\it Let} $f_1,f_2$ {\it satisfy} (8.47). {\it Then the convolution}
$$f_1*f_2\in {\cal C}(\R^3)\cap L^{\infty}(\R^3),\eqno(8.48)$$
{\it where}
$$(f_1*f_2)(p)=\int\limits_{\R^3}f_1(p-\xi)f_2(\xi)d\xi,\ \ p\in\R^3.
\eqno(8.49)$$

Lemma 8.1 follows from the following properties of (fixed)
$f_1\in L^{\infty}(\R^3)$, $f_2\in L^1(\R^3)$:
$$\eqalignno{
&\int\limits_{\R^3}|f_2(\xi)|d\xi<\infty,&(8.50a)\cr
&\int\limits_{{\cal B}_r}|f_2(\xi)|d\xi\to 0\ \ {\rm as}\ \ r\to +\infty
,&(8.50b)\cr
&\sup_{mes\,{\cal A}\le\ep}\int\limits_{{\cal A}}|f_2(\xi)|d\xi\to 0\ \
{\rm as}\ \ \ep\to 0,&(8.50c)\cr}$$
$$\eqalign{
&\forall\ r>0,\ep>0,\lambda>1\ \exists\ u\in {\cal C}({\cal B}_{r+1})\ \
{\rm such\ that}\cr
&mes\,supp\,(f_1-u)<\ep,\ \ \|u\|_{{\cal C}({\cal B}_{r+1})}\le\lambda
\|f_1\|_{L^{\infty}(\R^3)},\cr}\eqno(8.50d)$$
where
$${\cal B}_r=\{\xi\in\R^3:\ \ |\xi|<r\}.\eqno(8.51)$$
The proof of (3.2) is completed.

\vskip 2 mm
{\it Proof of (3.4).}
Due to (3.3a), we have that
$$\|(A(k)-A(l))U\|_{\mu}\le 2c_1(\mu)\|\hat v\|_{\mu}\|U\|_{\mu},\ \
k,l\in\Sigma.\eqno(8.52)$$
Besides, we have that
$$\eqalign{
&|(A(k)-A(l))U(p)|\le (\Delta_1(l,\ep,p)+\Delta_2(k,l,\ep,p)+
\Delta_3(k,l,\ep,r,p)+\Delta_4(k,l,r,p))\times\cr
&\|\hat v\|_{\mu}\|U\|_{\mu},\cr}\eqno(8.53)$$
where
$$\eqalignno{
&\Delta_1(l,\ep,p)=\int\limits_{{\cal D}(l,\ep)}{d\xi\over
(1+|p+\xi|)^{\mu}(1+|\xi|)^{\mu}|\xi^2+2l\xi|},&(8.54)\cr
&\Delta_2(k,l,\ep,p)=\int\limits_{{\cal D}(l,\ep)}{d\xi\over
(1+|p+\xi|)^{\mu}(1+|\xi|)^{\mu}|\xi^2+2k\xi|},&(8.55)\cr
&\Delta_3(k,l,\ep,r,p)=\int\limits_{{\cal B}_r\b {\cal D}(l,\ep)}
{2|(k-l)\xi|d\xi\over
(1+|p+\xi|)^{\mu}(1+|\xi|)^{\mu}|\xi^2+2k\xi||\xi^2+2l\xi|},&(8.56)\cr
&\Delta_4(k,l,r,p)=\int\limits_{\R^3\b {\cal B}_r}
{2|(k-l)\xi|d\xi\over
(1+|p+\xi|)^{\mu}(1+|\xi|)^{\mu}|\xi^2+2k\xi||\xi^2+2l\xi|},&(8.57)\cr}$$
where
$${\cal D}(l,\ep)=\{\xi\in\R^3:\ \ |\xi^2+2l\xi|\le\ep\},\eqno(8.58)$$
${\cal B}_r$ is defined by (8.51),
$$0<\ep\le 1,\ 2|l|+2\le r,\ |k-l|\le 1,\ k,l\in\Sigma,\ p\in\R^3,
\eqno(8.59)$$
where $|z|=(|Re\,z|^2+|Im\,z|^2)^{1/2}$ for $z\in\C^d$. Note that
$$\eqalign{
&|\xi^2+2l\xi|\ge |\xi^2+2Re\,l \xi|\ge |\xi|(|\xi|-2|Re\,l|)\ge 2|\xi|\ge
4>\ep\cr
&{\rm for}\ \ \xi\in\R^3\b {\cal B}_r\cr}\eqno(8.60)$$
and, therefore,
$${\cal D}(l,\ep)\subset {\cal B}_r \eqno(8.61)$$
under conditions (8.59). Further, we estimate separately $\Delta_1$,
$\Delta_2$, $\Delta_3$ and  $\Delta_4$.

\vskip 2 mm
{\it Estimate of} $\Delta_1$.
In a similar way with (8.6), (8.7) we obtain that
$$\Delta_1(l,\ep,p)\le \bigl(1+|p|/2\bigr)^{-\mu}\bigl(\Delta_{1,1}(l,\ep)+
(\Delta_{1,2}(l,\ep,p)\bigr),\eqno(8.62)$$
$$\eqalign{
&\Delta_{1,1}(l,\ep)=\int\limits_{{\cal D}(l,\ep)}{d\xi\over
(1+|\xi|)^{\mu}|\xi^2+2l\xi|},\cr
&\Delta_{1,2}(l,\ep,p)=\int\limits_{{\cal D}(l,\ep)}{d\xi\over
(1+|p+\xi|)^{\mu}|\xi^2+2l\xi|},\cr}\eqno(8.63)$$
where $0<\ep\le 1$, $l\in\Sigma$, $p\in\R^3$. In addition,
$$\Delta_{1,1}(l,\ep)\le\Delta_{1,3}(l,\ep),\ \
\Delta_{1,2}(l,\ep,p)\le\Delta_{1,3}(l,\ep),\eqno(8.64)$$
where
$$\eqalign{
&\Delta_{1,3}(l,\ep)=\int\limits_{{\cal D}(l,\ep)}{d\xi\over |\xi^2+2l\xi|}\le
\int\limits_{{\cal D}(l,\ep)}{\sqrt{2}d\xi\over
{|\xi^2+2Re\,l \xi|+2|Im\,l \xi|}}\buildrel (8.13) \over =\cr
&\int\limits_{|(\xi+Re\,l)^2-(Re\,l)^2+2i Im\,l(\xi+Re\,l)|\le\ep}
{\sqrt{2}d\xi\over {|(\xi+Re\,l)^2-(Re\,l)^2|+2|Im\,l(\xi+Re\,l)|}}\le\cr
&\int\limits_{|\xi^2-(Re\,l)^2|\le\ep}
{\sqrt{2}d\xi\over {|\xi^2-(Re\,l)^2|+2|Im\,l \xi|}}\buildrel
\rho_l=|Re\,l|,(8.13)\over \le
\int\limits_{|\xi^2-\rho_l^2|\le\ep}
{\sqrt{2}d\xi\over {|\xi^2-\rho_l^2|+2\rho_l|\xi_3|}}\le\cr
&\int\limits_{(\max\,(\rho_l^2-\ep,0))^{1/2}}^{(\rho_l^2+\ep)^{1/2}}
\int\limits_{-\pi}^{\pi}\int\limits_0^{\pi}
{\sqrt{2}r^2\sin\psi\,d\psi\,d\v\,dr\over
{|r^2-\rho_l^2|+2\rho_lr|\cos\psi|}}=4\pi\sqrt{2}\Delta_{1,4}(\rho_l,\ep),\cr}
\eqno(8.65)$$
$$\eqalign{
&\Delta_{1,4}(\rho,\ep)=
\int\limits_{(\max\,(\rho^2-\ep,0))^{1/2}}^{(\rho^2+\ep)^{1/2}}
\int\limits_0^{\pi/2}{r^2\sin\psi\,d\psi\,dr\over
{|r^2-\rho^2|+2\rho\,r\cos\psi}}\buildrel (8.29) \over = \cr
&{1\over 2}\int\limits_{(\max\,(\rho^2-\ep,0))^{1/2}}^{(\rho^2+\ep)^{1/2}}
\ln\bigl(1+{2\rho\,r\over
|r^2-\rho^2|}\bigr){r\over \rho}dr\le\cr
&{1\over 2}\biggl(\int\limits_{(\max\,(\rho^2-\ep,0))^{1/2}}^{\rho}+
\int\limits_{\rho}^{(\rho^2+\ep)^{1/2}}\biggr)
\ln\bigl(1+{2\rho \over |r-\rho|}\bigr)\bigl(1+{{r-\rho}\over \rho}\bigr)dr=
\cr
&\Delta_{1,4,1}(\rho,\ep)+\Delta_{1,4,2}(\rho,\ep),\cr}\eqno(8.66)$$
where  $\Delta_{1,4,1}$, $\Delta_{1,4,2}$ correspond to
$\int\limits_{(\max\,(\rho^2-\ep,0))^{1/2}}^{\rho}$,
$\int\limits_{\rho}^{(\rho^2+\ep)^{1/2}}$ respectively, $0<\ep\le 1$,
$l\in\Sigma$, $\rho_l=|Re\,l|=|l|/\sqrt{2}$, $0\le\rho$. In addition,
$$\Delta_{1,4}(\rho,\ep)=\int\limits_0^{\ep^{1/2}}dr=\ep^{1/2}\ \ {\rm for}
\ \ \rho=0,\ 0<\ep\le 1,\eqno(8.67a)$$
$$\eqalign{
&\Delta_{1,4,1}(\rho,\ep)\le {1\over 2}
\int\limits_{(\max\,(\rho^2-\ep,0))^{1/2}}^{\rho}
\ln\bigl(1+{2\rho \over |r-\rho|}\bigr)dr\buildrel (8.68) \over \le {1\over 2}
\int\limits_{(\max\,(\rho^2-\ep,0))^{1/2}}^{\rho}{1\over \alpha}
\bigl({2\rho\over {\rho-r}}\bigr)^{\alpha}dr=\cr
&-{(2\rho)^{\alpha}\over 2\alpha(1-\alpha)}
(\rho-r)^{1-\alpha}\bigg|_{(\max\,(\rho^2-\ep,0))^{1/2}}^{\rho}=\cr
&{(2\rho)^{\alpha}\over 2\alpha(1-\alpha)}
\bigl(\rho-(\max\,(\rho^2-\ep,0))^{1/2}\bigr)^{1-\alpha}\buildrel (8.69) \over
 \le {(2\rho)^{\alpha}\over 2\alpha(1-\alpha)}\ep^{(1-\alpha)/2}\cr}
\eqno(8.67b)$$
for $\rho>0$, $0<\ep\le 1$, $0<\alpha<1$,
$$\eqalign{
&\Delta_{1,4,2}(\rho,\ep)\buildrel (8.68) \over\le
\int\limits_{\rho}^{(\rho^2+\ep)^{1/2}}{1\over 2\alpha}
\bigl({2\rho\over {\eta-\rho}}\bigr)^{\alpha}dr+
\int\limits_{\rho}^{(\rho^2+\ep)^{1/2}}dr=\cr
&{(2\rho)^{\alpha}\over 2\alpha(1-\alpha)}
\bigl((\rho^2+\ep)^{1/2}-\rho\bigr)^{1-\alpha}+
\bigl((\rho^2+\ep)^{1/2}-\rho\bigr)\buildrel (8.69) \over \le\cr
&\biggl({(2\rho)^{\alpha}\over 2\alpha(1-\alpha)}+1\biggr)\ep^{(1-\alpha)/2}
\cr}\eqno(8.67c)$$
for $\rho>0$, $0<\ep\le 1$, $0<\alpha< 1$.

Note that in (8.67b), (8.67c) we used the inequalities
$$\eqalignno{
&\ln\,(1+x)\le\alpha^{-1}x^{\alpha}\ \ {\rm for}\ \ x\ge 0,\ \ 0<\alpha\le 1,
&(8.68)\cr
&(x+\ep)^{1/2}-x^{1/2}\le\ep^{1/2}\ \ {\rm for}\ \ x\ge 0,\ \ \ep\ge 0.
&(8.69)\cr}$$
Due to (8.62)-(8.67) we have that
$$\Delta_1(l,\ep,p)\le 8\pi\sqrt{2}(1+|p|/2)^{-\mu}
\biggl(1+{(2|Re\,l|)^{\alpha}\over \alpha(1-\alpha)}\biggr)\ep^{(1-\alpha)/2}
\eqno(8.70)$$
for $l\in\Sigma$, $0<\ep\le 1$, $p\in\R^3$, $0<\alpha<1$.

\vskip 2 mm
{\it Estimate of} $\Delta_2$.
In a similar way with (8.62)-(8.65) we obtain that
$$\Delta_2(k,l,\ep,p)\le 2(1+|p|/2)^{-\mu}\tilde\Delta_2(k,l,\ep),\eqno(8.71)
$$
$$\eqalign{
&\tilde\Delta_2(k,l,\ep)=\int\limits_{{\cal D}(l,\ep)}{\sqrt{2}d\xi\over
{|\xi^2+2Re\,k\xi|+2|Im\,k\xi|}}\le\cr
&\int\limits_{|((\xi+Re\,k)-(Re\,k-Re\,l))^2-(Re\,l)^2|\le\ep}{\sqrt{2}d\xi
\over {|(\xi+Re\,k)^2-(Re\,k)^2|+2|Im\,k(\xi+Re\,k)|}}=\cr
&\int\limits_{|(\xi-(Re\,k-Re\,l))^2-(Re\,l)^2|\le\ep}{\sqrt{2}d\xi
\over {|\xi^2-(Re\,k)^2|+2|Im\,k\xi|}},\cr}\eqno(8.72)$$
where $k,l\in\Sigma$, $0<\ep\le 1$, $p\in\R^3$. Note that
$$\eqalignno{
&|(\xi-\zeta)^2-\rho^2|\le\ep\Leftrightarrow \max\,(\rho^2-\ep,0)\le
(\xi-\zeta)^2\le\rho^2+\ep\Rightarrow &(8.73)\cr
&\max\,((\max\,(\rho^2-\ep,0))^{1/2}-|\zeta|,0)\le |\xi|\le |\zeta|+
(\rho^2+\ep)^{1/2} &(8.74)\cr
&\buildrel (8.69) \over \Rightarrow \max\,(\rho-\ep^{1/2}-|\zeta|,0)\le
|\xi|\le \rho+\ep^{1/2}+|\zeta|,&(8.75)\cr}$$
where $\xi,\zeta\in\R^3$, $\rho\ge 0$, $0<\ep\le 1$. Using (8.72) and
(8.73)-(8.75) for $\zeta=Re\,k-Re\,l$, $\rho=|Re\,l|$,
 in a similar way with (8.65), (8.66) we obtain that
$$\eqalign{
&\tilde\Delta_2(k,l,\ep)=\cr
&2\pi\sqrt{2}\biggl(\int\limits_{\max\,(\rho_l-\delta,0)}^{\rho_k}+
\int\limits_{\rho_k}^{\rho_l+\delta}\biggr)
\ln\,\biggl(1+{2\rho_k\over |r-\rho_k|}\biggr)\biggl(1+{{r-\rho_k}\over
\rho_k}\biggr)dr=\cr
&2\pi\sqrt{2}(\tilde\Delta_{2,1}(\rho_k,\rho_l,\delta)+
\tilde\Delta_{2,2}(\rho_k,\rho_l,\delta)),\cr
&\rho_k=|Re\,k|\ne 0,\ \ \rho_l=|Re\,l|,\ \ \delta=\ep^{1/2}+|Re\,k-Re\,l|,\ \
 k,l\in\Sigma,\ \ 0<\ep\le 1,\cr}\eqno(8.76)$$
where  $\tilde\Delta_{2,1}$, $\tilde\Delta_{2,2}$  correspond to
$\int\limits_{\max\,(\rho_l-\delta,0)}^{\rho_k}$,
$\int\limits_{\rho_k}^{\rho_l+\delta}$ respectively. In addition, in a
similar way with (8.67) we obtain that
$$\tilde\Delta_2(k,l,\ep)\le 4\pi\sqrt{2}(\delta+\rho_l-\rho_k)\le
4\pi\sqrt{2}(\ep^{1/2}+2|Re\,k-Re\,l|)\ \ {\rm for}\ \ k=0,\eqno(8.77a)$$
$$\eqalign{
&\tilde\Delta_{2,1}(\rho_k,\rho_l,\delta)\le {(2\rho_k)^{\alpha}\over
\alpha(1-\alpha)}(\rho_k-\max\,(\rho_l-\delta,0))^{1-\alpha}\le\cr
&{(2\rho_k)^{\alpha}\over \alpha(1-\alpha)}(\delta+|\rho_k-\rho_l|)^{1-\alpha}
\le
{(2|Re\,k|)^{\alpha}\over \alpha(1-\alpha)}
(\ep^{1/2}+2|Re\,k-Re\,l|)^{1-\alpha},\cr}\eqno(8.77b)$$
$$\eqalign{
&\tilde\Delta_{2,2}(\rho_k,\rho_l,\delta)\le {(2\rho_k)^{\alpha}\over
\alpha(1-\alpha)}(\rho_l+\delta-\rho_k)^{1-\alpha}+
2(\rho_l+\delta-\rho_k)\le\cr
&{(2|Re\,k|)^{\alpha}\over \alpha(1-\alpha)}
(\ep^{1/2}+2|Re\,k-Re\,l|)^{1-\alpha}+2(\ep^{1/2}+2|Re\,k-Re\,l|),\cr}
\eqno(8.77c)$$
where $k,l,\ep,\rho_k,\rho_l,\delta$ are the same as in (8.76) and
$0<\alpha<1$.

Due to (8.71), (8.76), (8.77) we have that
$$\eqalign{
&\Delta_2(k,l,\ep,p)\le 8\pi\sqrt{2}(1+|p|/2)^{-\mu}
\biggl({(2|Re\,k|)^{\alpha}\over \alpha(1-\alpha)}+3^{\alpha}\biggr)
(\ep^{1/2}+2|Re\,k-Re\,l|)^{1-\alpha}\cr
&{\rm for}\ \ k,l\in\Sigma,\ \ |k-l|\le 1,\ \ 0<\ep\le 1,\ \ p\in\R^3,\ \
0<\alpha<1.\cr}\eqno(8.78)$$

\vskip 2 mm
{\it Estimate of} $\Delta_3$.
We have that
$$\eqalign{
&\Delta_3(k,l,\ep,r,p)\buildrel (8.51), (8.58) \over\le
\int\limits_{{\cal B}_r}{2|k-l|rd\xi\over
(1+|p+\xi|)^{\mu}(1+|\xi|)^{\mu}|\xi^2+2k\xi|\ep}\buildrel (8.3a)\over\le\cr
&{2|k-l|r\over \ep}{c_1(\mu)\over (1+|p|)^{\mu}}\cr}\eqno(8.79)$$
under conditions (8.59).

\vskip 2 mm
{\it Estimate of} $\Delta_4$.
We have that
$$\eqalign{
&\Delta_4(k,l,r,p)\buildrel (8.57), (8.59), (8.60) \over\le
\int\limits_{\R^3\b {\cal B}_r}{|k-l|d\xi\over
(1+|p+\xi|)^{\mu}(1+|\xi|)^{\mu}|\xi^2+2k\xi|}\buildrel (8.3a)\over\le\cr
&|k-l|{c_1(\mu)\over (1+|p|)^{\mu}}\cr}\eqno(8.80)$$
under conditions (8.59).

Now formulas (3.4) follow from (8.52), (8.53) and estimates (8.70),
(8.78)-(8.80) with $\ep=|k-l|^{\beta}$, $0<|k-l|\le 1$ for fixed
$k\in\Sigma$, $r\ge 2(|k|+\sqrt{2})+2$, $\alpha\in ] 0,1 [$ and
$\beta\in ] 0,1 [$.

The proof of (3.4) is completed.

Finally, property (3.5) follows from the presentation
$$\eqalign{
&(A(k)U)(p)-(A(k^{\prime})U)(p^{\prime})=\cr
&((A(k)U)(p)-(A(k)U)(p^{\prime}))+((A(k)U)(p^{\prime})-
(A(k^{\prime})U)(p^{\prime}))\cr}\eqno(8.81)$$
and properties (3.2), (3.4).
The proof of Lemma 3.1 is completed.

\vskip 2 mm
{\it Proof of Proposition 3.1.}
Proposition 3.1 follows from equation (1.5) written as
$$H(k,\cdot)=\hat v-A(k)H(k,\cdot)\ \eqno(8.82)$$
and Lemma 3.1. In addition, to obtain (3.8a), (3.12a) we use the presentation
$$\tilde H(k,p)-\tilde H(k^{\prime},p^{\prime})=
(\tilde H(k,p)-\tilde H(k,p^{\prime}))+
(\tilde H(k,p^{\prime})-\tilde H(k^{\prime},p^{\prime})),\eqno(8.83)$$
where
$$\eqalignno{
&\tilde H(k,\cdot)\buildrel \rm def \over = H(k,\cdot)-\hat v\buildrel (8.82)
\over = -A(k)H(k,\cdot),&(8.84)\cr
&\tilde H(k,\cdot)\buildrel (3.2),(8.84) \over \in C(\R^3)\ \ {\rm as\ soon\
as}\ \ H(k,\cdot)\in L^{\infty}_{\mu}(\R^3),&(8.85)\cr}$$
$$\eqalign{
&\tilde H(k,\cdot)-\tilde H(k^{\prime},\cdot)=H(k,\cdot)-H(k^{\prime},\cdot)
\buildrel (8.82) \over = ((I+A(k))^{-1}-(I+A(k^{\prime}))^{-1})\hat v=\cr
&(I+A(k))^{-1}((I+A(k^{\prime}))^{-1}-(I+A(k)))
(I+A(k^{\prime}))^{-1})\hat v=\cr
&(I+A(k))^{-1}(A(k^{\prime})-A(k))(I+A(k^{\prime}))^{-1})\hat v,\cr}
\eqno(8.86)$$
$$\tilde H(k,\cdot)-\tilde H(k^{\prime},\cdot)\buildrel (8.85) \over \in
C(\R^3)\ \ {\rm as\ soon\ as}\ \ H(k,\cdot),H(k^{\prime},\cdot)
\in L^{\infty}_{\mu}(\R^3),\eqno(8.87a)$$
$$\eqalign{
&\|\tilde H(k,\cdot)-\tilde H(k^{\prime},\cdot)\|_{\mu}\buildrel
(3.4),(8.86) \over \to 0\ \ {\rm as}\ \ k^{\prime}\to k\cr
&{\rm as\ soon\ as}\ \ (I+A(k^{\prime}))^{-1}\ \ {\rm is\ uniformly\ bounded\
 in\ a\ neighborhood\ of}\ \ k,\cr}\eqno(8.87b)$$
$$\eqalign{
&\sup_{p^{\prime}\in\R^3}(1+|p^{\prime}|)^{\mu}
|\tilde H(k,p^{\prime})-\tilde H(k^{\prime},p^{\prime})|\buildrel (8.87)
\over \to 0\ \ {\rm as}\ \ k^{\prime}\to k\cr
&{\rm as\ soon\ as}\ \ (I+A(k^{\prime}))^{-1}\ \ {\rm is\ uniformly\ bounded\
in\ a\ neighborhood\ of}\ \ k,\cr}\eqno(8.88)$$
where $k$, $k^{\prime}\in\Sigma$, $p,p^{\prime}\in\R^3$.

The proof of Proposition 3.1 is completed.

\vskip 2 mm
{\bf 9. Proof of Lemma 5.1}

The proof of Lemma 5.1 of the present work is similar to the proof of
Lemma 4.1 of [No5]. Proceeding from (3.13), (4.3), (4.4a), (4.7), (5.2), (5.3)
 in a similar way with the  proof of Lemma 4.1 of [No5] we obtain that:
$$\eqalign{
&{\pa\over \pa\bar\lambda}H(k(\lambda,p),p)=-{\pi\over 2}
\int\limits_{\{\xi\in\R^3:\ \xi^2+2k\xi=0\}}
\biggl({\pa\bar\kappa_1\over \pa\bar\lambda}\theta\xi+
{\pa\bar\kappa_2\over \pa\bar\lambda}\omega\xi\biggr)\times\cr
&H(k,-\xi)H(k+\xi,p+\xi){ds\over |Im\,k|^2},\ \
\lambda\in\C\b 0,\ \ p\in\R^3\b {\cal L}_{\nu},\cr}\eqno(9.1)$$
where $k=k(\lambda,p)$, $\kappa_1=\kappa_1(\lambda,p)$,
$\kappa_2=\kappa_2(\lambda,p)$ are defined in (4.7), $\theta=\theta(p)$,
$\omega=\omega(p)$ are the vector-functions of (4.3), (4.4a), $ds$ is
arc-length measure on the circle   $\{\xi\in\R^3:\ \xi^2+2k\xi=0\}$ and, in
addition,
$$ds=|Re\,k|d\v, \eqno(9.2)$$
$$\eqalign{
&{\pa\bar\kappa_1\over \pa\bar\lambda}\theta\xi+
{\pa\bar\kappa_2\over \pa\bar\lambda}\omega\xi=
\biggl({\pa\bar\kappa_1\over \pa\bar\lambda}Re\,\kappa_1+
{\pa\bar\kappa_2\over \pa\bar\lambda}Re\,\kappa_2\biggr)(\cos\v-1)+\cr
&{|p|\over 2|Im\,k|}
\biggl({\pa\bar\kappa_1\over \pa\bar\lambda}Im\,\kappa_2-
{\pa\bar\kappa_2\over \pa\bar\lambda}Im\,\kappa_1\biggr)\sin\v\cr}\eqno(9.3)$$
under the assumption that the circle   $\{\xi\in\R^3:\ \xi^2+2k\xi=0\}$  is
parametrized by $\v\in ] -\pi,\pi [$ according to (5.2). (Note that in the
proof of Lemma 4.1 of [No5] the $\bar\pa$-equation similar to (9.1) is not
valid for $|\lambda|=1$ but it is not indicated because of a misprint.)

The difinition of $\kappa_1$, $\kappa_2$ (see (4.7)) implies that
$${\pa\bar\kappa_1\over \pa\bar\lambda}=
-{i|p|\over 4}\bigl(1-{1\over \bar\lambda^2}\bigr),\ \
{\pa\bar\kappa_2\over \pa\bar\lambda}=
{|p|\over 4}\bigl(1+{1\over \bar\lambda^2}\bigr),\eqno(9.4)$$
$$\eqalign{
&Re\,\kappa_1={i|p|\over 8}\bigl(\lambda-\bar\lambda+{1\over \lambda}-
{1\over \bar\lambda}\bigr),\ \
Im\,\kappa_1={|p|\over 8}\bigl(\lambda+\bar\lambda+{1\over \lambda}+
{1\over \bar\lambda}\bigr),\cr
&Re\,\kappa_2={|p|\over 8}\bigl(\lambda+\bar\lambda-{1\over \lambda}-
{1\over \bar\lambda}\bigr),\ \
Im\,\kappa_2={|p|\over 8i}\bigl(\lambda-\bar\lambda-{1\over \lambda}+
{1\over \bar\lambda}\bigr),\cr}\eqno(9.5)$$
where $\lambda\in\C\b 0$, $p\in\R^3$. Due to (9.4), (9.5) we have that
$$\eqalign{
&{\pa\bar\kappa_1\over \pa\bar\lambda}Re\,\kappa_1+
{\pa\bar\kappa_2\over \pa\bar\lambda}Re\,\kappa_2=\cr
&{|p|^2\over 32}\bigl(\bigl(1-{1\over \bar\lambda^2}\bigr)
\bigl(\lambda-\bar\lambda+{1\over \lambda}-{1\over \bar\lambda}+
\bigl(1+{1\over \bar\lambda^2}\bigr)
\bigl(\lambda+\bar\lambda-{1\over \lambda}-{1\over \bar\lambda}\bigr)\bigr)=
\cr
&{|p|^2\over 32}\bigl(\lambda-\bar\lambda+{1\over \lambda}-
{1\over \bar\lambda}-{\lambda\over \bar\lambda^2}+{1\over \bar\lambda}-
{1\over \lambda\bar\lambda^2}+{1\over \bar\lambda^3}+\cr
&\lambda+\bar\lambda-{1\over \lambda}-
{1\over \bar\lambda}+{\lambda\over \bar\lambda^2}+{1\over \bar\lambda}-
{1\over \lambda\bar\lambda^2}-{1\over \bar\lambda^3}\bigr)=\cr
&{|p|^2\over 16}\bigl(\lambda-{1\over \lambda\bar\lambda^2}\bigr)=
{|p|^2\over 16}\lambda\bigl(1-{1\over |\lambda|^4}\bigr),\cr}\eqno(9.6)$$
$$\eqalign{
&{\pa\bar\kappa_1\over \pa\bar\lambda}Im\,\kappa_2-
{\pa\bar\kappa_2\over \pa\bar\lambda}Im\,\kappa_1=\cr
&-{|p|^2\over 32}\bigl(\bigl(1-{1\over \bar\lambda^2}\bigr)
\bigl(\lambda-\bar\lambda-{1\over \lambda}+{1\over \bar\lambda}+
\bigl(1+{1\over \bar\lambda^2}\bigr)
\bigl(\lambda+\bar\lambda+{1\over \lambda}+{1\over \bar\lambda}\bigr)\bigr)=
\cr
&-{|p|^2\over 32}\bigl(\lambda-\bar\lambda-{1\over \lambda}+
{1\over \bar\lambda}-{\lambda\over \bar\lambda^2}+{1\over \bar\lambda}+
{1\over \lambda\bar\lambda^2}-{1\over \bar\lambda^3}+\cr
&\lambda+\bar\lambda+{1\over \lambda}+
{1\over \bar\lambda}+{\lambda\over \bar\lambda^2}+{1\over \bar\lambda}+
{1\over \lambda\bar\lambda^2}+{1\over \bar\lambda^3}\bigr)=\cr
&-{|p|^2\over 16}\bigl(\lambda+{2\over \bar\lambda}+
{1\over \lambda\bar\lambda^2}\bigr)=
-{|p|^2\over 16}{(|\lambda|^2+1)^2\over |\lambda|^2\bar\lambda}.\cr}
\eqno(9.7)$$
Due to (9.6), (9.7), (4.8) we have that
$$\eqalignno{
&\biggl({\pa\bar\kappa_1\over \pa\bar\lambda}Re\,\kappa_1+
{\pa\bar\kappa_2\over \pa\bar\lambda}Re\,\kappa_2\biggr){1\over |Im\,k|}=
{|p|\over 4}{(|\lambda|^2-1)\over \bar\lambda|\lambda|},&(9.8)\cr
&\biggl({\pa\bar\kappa_1\over \pa\bar\lambda}Im\,\kappa_2-
{\pa\bar\kappa_2\over \pa\bar\lambda}Im\,\kappa_1\biggr)
{|p|\over 2|Im\,k|^2}=-{|p|\over 2\bar\lambda},&(9.9)\cr}$$
where $(\lambda,p)\in (\C\b 0)\times (\R^3\b {\cal L}_{\nu})$.

The $\bar\pa$-equation (5.1) follows from (9.1), (9.2), (9.3), (9.8), (9.9)
and the property that $|Re\,k|=|Im\,k|$ for $k\in\Sigma$ defined by (1.7).

Lemma 5.1 is proved.

\vskip 2 mm
{\bf 10. Proof of Lemma 5.2}

Let us show, first, that
$$\{U_1,U_2\}\in  L_{local}^{\infty}((\C\b 0)\times (\R^3\b {\cal L}_{\nu})).
\eqno(10.1)$$
Property (10.1) follows from definition (5.5), the properties
$$\eqalign{
&U_1(k,-\xi(k,\v))\in L^{\infty}(\Sigma\times [0,2\pi])\ \ {\rm (as\ a\
function\ of}\ \ k,\v),\cr
&U_1(k,-\xi(k,\v))\in L^{\infty}(\Omega\times [0,2\pi])\ \ {\rm (as\ a\
function\ of}\ \ k,p,\v\cr
&({\rm with\ no\ dependence\ on}\ \ p)),\cr}\eqno(10.2)$$
$$U_2(k+\xi(k,\v),p+\xi(k,\v))\in L^{\infty}(\Omega\times [0,2\pi])\ \
{\rm (as\ a\ function\ of}\ \ k,p,\v),\eqno(10.3)$$
where
$$\eqalignno{
&\Sigma=\{k\in\C^3:\ k^2=0\},\ \
\Omega=\{k\in\C^3,\ p\in\R^3:\ k^2=0,\ p^2=2kp\},&(10.4)\cr
&\xi(k,\v)=Re\,k(\cos\v-1)+k^{\perp}\sin\v,\ \
k^{\perp}={Im\,k\times Re\,k\over |Im\,k|}&(10.5)\cr}$$
(where $\times$ in (10.5) denotes vector product), and from Lemma 4.1. In
turn,
(10.2) follows from $U_1\in L^{\infty}(\Omega)$, definition (10.4) and the
fact
that $p=-\xi(k,\v)$, $\v\in [0,2\pi]$, is a parametrization of the set
$\{p\in\R^3:\ p^2=2kp\}$, $k\in\Sigma\b \{0\}$. To prove (10.3), consider
$$\Theta=\{k\in\C^3,\ l\in\C^3:\ k^2=l^2=0,\ \ Im\,k=Im\,l\}.\eqno(10.6)$$
Note that
$$\eqalign{
&\Theta\approx\Omega,\cr
&(k,l)\in\Theta\Rightarrow (k,k-l)\in\Omega,\ \
(k,p)\in\Omega\Rightarrow (k,k-p)\in\Theta.\cr}\eqno(10.7)$$
Consider
$$u_2(k,l)=U_2(k,k-l),\ \ (k,l)\in\Theta.\eqno(10.8)$$
The property $U_2\in L^{\infty}(\Omega)$ is equivalent to the property
$u_2\in L^{\infty}(\Theta)$. Property (10.3) is equivalent to the property
$$u_2(k+\xi(k,\v),l)\in L^{\infty}(\Theta\times [0,2\pi])\ \ {\rm (as\ a\
function\ of}\ \ k,l,\v).\eqno(10.9)$$
Property (10.9) follows from the property
$$\eqalign{
&u_2(\zeta(l,\psi,\v)+iIm\,l,l)\in L^{\infty}(\Sigma\times [0,2\pi]\times
 [0,2\pi])\cr
&{\rm (as\ a\ function\ of}\ \ l,\psi,\v),\cr}\eqno(10.10)$$
where
$$\zeta(l,\psi,\v)=Re\,l\cos(\v-\psi)+l^{\perp}\sin(\v-\psi),
\ \ l^{\perp}={Im\,l\times Re\,l\over |Im\,l|} \eqno(10.11)$$
(where $\times$ in (10.11) denotes vector product). Note that
$k=\zeta(l,\psi,\v)$, $\v\in [0,2\pi]$ at fixed $\psi\in [0,2\pi]$ is a
parametrization of the set
$S_l=\{k\in\C^3:\ k^2=l^2,\ Im\,k=Im\,l\}$, $l\in\Sigma\b 0$.
In turn, (10.10) follows from $u_2\in L^{\infty}(\Theta)$, definition (10.6)
and
the aforementioned fact concerning the parametrization of $S_l$. Thus,
properties (10.10), (10.9), (10.3) are proved. This completes the proof of
(10.1).

Let us prove now (5.8).

We have that
$$\{U_1,U_2\}=\{U_1,U_2\}_1+\{U_1,U_2\}_2,\eqno(10.12)$$
where
$$\{U_1,U_2\}_1(\lambda,p)=-{\pi |p|(|\lambda|^2-1)\over
8\bar\lambda|\lambda|}\{U_1,U_2\}_3(\lambda,p),\eqno(10.13a)$$
$$\eqalign{
&\{U_1,U_2\}_3(\lambda,p)=\int_{-\pi}^{\pi}(\cos\v-1)\times\cr
&U_1(k(\lambda,p),-\xi(\lambda,p,\v))U_2(k(\lambda,p)+\xi(\lambda,p,\v),
p+\xi(\lambda,p,\v))d\v,\cr}\eqno(10.13b)$$
$$\{U_1,U_2\}_2={\pi |p|\over 4\bar\lambda}
\{U_1,U_2\}_4(\lambda,p),\eqno(10.14a)$$
$$\eqalign{
&\{U_1,U_2\}_4(\lambda,p)=\int_{-\pi}^{\pi}\sin\v\times\cr
&U_1(k(\lambda,p),-\xi(\lambda,p,\v))U_2(k(\lambda,p)+\xi(\lambda,p,\v),
p+\xi(\lambda,p,\v))d\v,\cr}\eqno(10.14b)$$
$\lambda\in\C\b 0$, $p\in\R^3\b {\cal L}_{\nu}$.

Formulas (5.2), (5.3) imply that
$$|\xi|^2=|Re\,k|^2((\cos\v-1)^2+(\sin\v)^2)=4|Re\,k|^2(\sin\,(\v/2))^2,
\eqno(10.15)$$
where $\xi=\xi(\lambda,p,\v)$, $k=k(\lambda,p)$.

The relation $p^2=2k(\lambda,p)p$, $\lambda\in\C\b 0$, $p\in\R^3\b
{\cal L}_{\nu}$, implies that
$$p=-Re\,k(\lambda,p)(\cos\psi-1)-k^{\perp}(\lambda,p)\sin\psi \eqno(10.16)$$
for some $\psi=\psi(\lambda,p)\in [-\pi,\pi]$, where $k^{\perp}(\lambda,p)$
is defined by (5.3). Formulas (5.2), (5.3), (10.16) imply that
$$\eqalign{
&|p+\xi|^2=|Re\,k|^2((\cos\v-\cos\psi)^2+(\sin\v-\sin\psi)^2)=
4|Re\,k|^2\bigl(\sin{{\v-\psi}\over 2}\bigr)^2,\cr
&|p|^2=4|Re\,k|^2\bigl(\sin{\psi\over 2}\bigr)^2,\cr}\eqno(10.17)$$
where $\xi=\xi(\lambda,p,\v)$, $k=k(\lambda,p)$, $\psi=\psi(\lambda,p)$.

Using the assumptions of Lemma 5.2 and formulas (10.13b), (10.14b), (10.15),
(10.17) we obtain that
$$\eqalign{
&|\{U_1,U_2\}_3(\lambda,p)|\le A(r,\psi,\mu,\mu)|||U_1|||_{\mu}
|||U_2|||_{\mu},\cr
&|\{U_1,U_2\}_4(\lambda,p)|\le B(r,\psi,\mu,\mu)|||U_1|||_{\mu}
|||U_2|||_{\mu}\cr}\eqno(10.18)$$
for $r=|Re\,k(\lambda,p)|$, $\psi=\psi(\lambda,p)$ (of (10.16)) and almost all
 $(\lambda,p)\in (\C\b 0)\times (\R^3\b {\cal L}_{\nu})$,  where
$$\eqalignno{
&A(r,\psi,\alpha,\beta)=\int_{-\pi}^{\pi}{(1-\cos\v)d\v\over
(1+2r|\sin(\v/2)|)^{\alpha}(1+2r|\sin({{\v-\psi}\over 2})|)^{\beta}},&(10.19a)
\cr
&B(r,\psi,\alpha,\beta)=\int_{-\pi}^{\pi}{|\sin\v|d\v\over
(1+2r|\sin(\v/2)|)^{\alpha}(1+2r|\sin({{\v-\psi}\over 2})|)^{\beta}},&(10.19b)
\cr
}$$
for $r\ge 0$, $\psi\in [-\pi,\pi]$, $\alpha\ge 2$, $\beta\ge 2$.
In addition, in (10.18) we used also that, in view of Lemma 4.1, properties
(10.2), (10.3) and definitions (10.13), (10.14), the variations of
$U_1$, $U_2$ on the sets of zero measure in $\Omega$ imply variations of
$\{U_1,U_2\}_3$ and $\{U_1,U_2\}_4$ on sets of zero measure, only, in
$(\C\b 0)\times (\R^3\b {\cal L}_{\nu})$.

Further, we use the following lemma of [No5].

\vskip 2 mm
{\bf Lemma 10.1} ([No5]). {\it Let} $r\ge 0$, $\psi\in [-\pi,\pi]$,
$\rho=2r|\sin(\psi/2)|$, $\alpha\ge 2$, $\beta\ge 2$. {\it Then}
$$\eqalignno{
&A(r,\psi,\alpha,\beta)\le\sum_{j=1}^4A_j(r,\psi,\alpha,\beta),&(10.20)\cr
&A_1(r,\psi,\alpha,\beta)\le\min\,\bigl({\rho^3\over 6r^3},{\rho\over r^3}
\bigr){1\over (1+\rho/2)^{\beta}},&(10.21)\cr
&A_2(r,\psi,\alpha,\beta)\le {\rho^3\over r^3}
{1\over (1+\rho/2)^{\alpha+1}},&(10.22)\cr
&A_3(r,\psi,\alpha,\beta)\le {4\rho^3\over r^3}
{1\over (1+\rho)^{\alpha}(1+\rho/2)},&(10.23)\cr
&A_4(r,\psi,\alpha,\beta)\le \bigl({3\over {1+r^2}}+
{2\pi\over (1+\sqrt{2}r)^{\alpha}}\bigr)
{1\over (1+\rho/2)^{\beta}},&(10.24)\cr
&B(r,\psi,\alpha,\beta)\le\sum_{j=1}^4B_j(r,\psi,\alpha,\beta),&(10.25)\cr
&B_1(r,\psi,\alpha,\beta)\le\min\,\bigl({\rho^2\over 2r^2},
{\sqrt{2}\rho\over r^2}\bigr){1\over (1+\rho/2)^{\beta}},&(10.26)\cr
&B_2(r,\psi,\alpha,\beta)\le {2\rho^2\over r^2}
{1\over (1+\rho/2)^{\alpha+1}},&(10.27)\cr
&B_3(r,\psi,\alpha,\beta)\le {4\rho^2\over r^2}
{1\over (1+\rho)^{\alpha}(1+\rho/2)},&(10.28)\cr
&B_4(r,\psi,\alpha,\beta)\le \bigl({5\over {1+r}}+
{3\over (1+\sqrt{2}r)^{\alpha}}\bigr)
{1\over (1+\rho/2)^{\beta}}.&(10.29)\cr}$$

\vskip 2 mm
{\bf Lemma 10.2.}
{\it Let}
$$r=r(\lambda,p)={\rho\over 4}\bigl(|\lambda|+{1\over |\lambda|}\bigr),\ \
|\sin\,(\psi/2)|={\rho\over 2r},\eqno(10.30)$$
{\it where} $\lambda\in\C\b 0$, $\rho\ge 0$, $\psi\in [-\pi,\pi]$. {\it Then}:
$$\eqalignno{
&{\rho\,||\lambda|^2-1|\over |\lambda|^2}A_1\le {4^3|\lambda|\over
\sqrt{6}(|\lambda|^2+1)^2(1+\rho/2)^{\beta}},&(10.31)\cr
&{\rho\,||\lambda|^2-1|\over |\lambda|^2}A_2\le {2\cdot 4^3|\lambda|\over
(|\lambda|^2+1)^2(1+\rho/2)^{\alpha}},&(10.32)\cr
&{\rho\,||\lambda|^2-1|\over |\lambda|^2}A_3\le {2\cdot 4^4|\lambda|\over
(|\lambda|^2+1)^2(1+\rho)^{\alpha}},&(10.33)\cr
&{\rho\,||\lambda|^2-1|\over |\lambda|^2}A_4\le {4\pi\rho||\lambda|^2-1|\over
|\lambda|^2(1+(\rho/4)(|\lambda|+|\lambda|^{-1}))^2
(1+\rho/2)^{\beta}},&(10.34)\cr
&{\rho\over |\lambda|}B_1\le {16\sqrt{2}|\lambda|\over
(|\lambda|^2+1)^2(1+\rho/2)^{\beta}},&(10.35)\cr
&{\rho\over |\lambda|}B_2\le {4^3|\lambda|\over
(|\lambda|^2+1)^2(1+\rho/2)^{\alpha}},&(10.36)\cr
&{\rho\over |\lambda|}B_3\le {2\cdot 4^3|\lambda|\over
(|\lambda|^2+1)^2(1+\rho)^{\alpha}},&(10.37)\cr
&{\rho\over |\lambda|}B_4\le {8\rho\over
|\lambda|(1+(\rho/4)(|\lambda|+|\lambda|^{-1}))
(1+\rho/2)^{\beta}},&(10.38)\cr}$$
{\it where} $A_j=A_j(r,|\psi|,\alpha,\beta)$,
$B_j=B_j(r,|\psi|,\alpha,\beta)$  {\it are the same as in Lemma} 10.1,
$j=1,2,3,4$, $\alpha\ge 2$, $\beta\ge 2$.

\vskip 2 mm
{\it Proof of Lemma 10.2.}
Using (10.30) we obtain that
$$\eqalignno{
&\rho\min\,\biggl({\rho^3\over 6r^3},{\rho\over r^3}\biggr)={4^3|\lambda|^3
\over (|\lambda|^2+1)^3}\min\,\biggl({\rho\over 6},{1\over \rho}\biggr)\le
{4^3|\lambda|^3\over \sqrt{6}(|\lambda|^2+1)^3},&(10.39)\cr
&\rho\min\,\biggl({\rho^2\over 2r^2},{\sqrt{2}\rho\over r^2}\biggr)=
{16|\lambda|^2
\over (|\lambda|^2+1)^2}\min\,\biggl({\rho\over 2},\sqrt{2}\biggr)\le
{16\sqrt{2}|\lambda|^2\over (|\lambda|^2+1)^2},&(10.40)\cr}$$
where $\lambda\in\C\b 0$, $\rho\ge 0$. Estimates (10.31), (10.35) follow from
(10.21), (10.26) and (10.39), (10.40). Estimates (10.32), (10.33), (10.36),
(10.37)
follow from (10.22), (10.23), (10.27), (10.28) and (10.30). Estimates (10.34),
 (10.38) follow from (10.24), (10.29), the inequalities
$$\eqalign{
&{3\over {1+r^2}}+{2\pi\over (1+\sqrt{2}r)^{\alpha}}\le
{4\pi\over (1+r)^2},\cr
&{5\over {1+r}}+{3\over (1+\sqrt{2}r)^{\alpha}}\le
{8\over {1+r}},\cr}\eqno(10.41)$$
where $r\ge 0$, $\alpha\ge 2$, and from (10.30). Lemma 10.2 is proved.

Estimate (5.8) follows from (10.12)-(10.14), (10.18), (4.8), (10.17) (for
$|p|$) and Lemmas 10.1, 10.2. Property (5.7) follows from (10.1) and (5.8).

Lemma 5.2 is proved.

\vskip 2 mm
{\it 11. Proof of Lemma 6.1}

Let
$$\eqalignno{
&J_1(\lambda)=\int_{\C}{|\zeta|\over (|\zeta|^2+1)^2}{d\,Re\zeta\,d\,Im\zeta
\over |\zeta-\lambda|},&(11.1)\cr
&J_2(\lambda,\rho)=\int_{\C}{(|\zeta|^2+1)\rho\over |\zeta|^2
(1+\rho(|\zeta|+|\zeta|^{-1}))^2}
{d\,Re\zeta\,d\,Im\zeta\over |\zeta-\lambda|},&(11.2)\cr
&J_3(\lambda,\rho)=\int_{\C}{\rho\over |\zeta|
(1+\rho(|\zeta|+|\zeta|^{-1}))}
{d\,Re\zeta\,d\,Im\zeta\over |\zeta-\lambda|},&(11.3)\cr}$$
where $\lambda\in\C$, $\rho>0$.

\vskip 2 mm
{\bf Lemma 11.1.}
{\it The following estimates hold}:
$$\eqalignno{
&J_1(\lambda)\le n_1,\ \ \lambda\in\C,&(11.4)\cr
&J_2(\lambda,\rho)\le n_2,\ \ \lambda\in\C,\ \ \rho>0,&(11.5)\cr
&J_3(\lambda,\rho)\le n_3,\ \ \lambda\in\C,\ \ \rho>0,&(11.6)\cr}$$
{\it for some positive constants} $n_1,n_2,n_3$ ({\it where} $J_1$, $J_2$,
$J_3$ {\it are defined by} (11.1)-(11.3)).

\vskip 2 mm
{\it Proof of Lemma 11.1}.

{\it Proof of (11.4)}. We have that
$$\eqalign{
&J_1(\lambda)\le \biggl(\int\limits_{|\zeta|\le |\zeta-\lambda|}+
\int\limits_{|\zeta|\ge |\zeta-\lambda|}\biggr){2|\zeta|\over
(|\zeta|^2+1)(|\zeta|+1)^2}
{d\,Re\zeta\,d\,Im\,\zeta\over |\zeta-\lambda|}\le\cr
&\int\limits_{|\zeta|\le |\zeta-\lambda|}
{2d\,Re\zeta\,d\,Im\,\zeta\over (|\zeta|^2+1)(|\zeta|+1)^2}+
\int\limits_{|\zeta|\ge |\zeta-\lambda|}
{2d\,Re\zeta\,d\,Im\,\zeta\over
(|\zeta-\lambda|^2+1)(|\zeta-\lambda|+1)|\zeta-\lambda|}\le\cr
&\int_0^{+\infty}{4\pi rdr\over (r^2+1)(r+1)^2}+
\int_0^{+\infty}{4\pi rdr\over (r^2+1)(r+1)r}\le n_1,\cr}\eqno(11.7)$$
where $\lambda\in\C$. Estimate (11.4) is proved.

{\it Proof of (11.5)}. We have that
$$J_2(\lambda,\rho)=J_{2,1}(\lambda,\rho)+J_{2,2}(\lambda,\rho),\eqno(11.8a)$$
$$\eqalign{
&J_{2,1}(\lambda,\rho)=\int\limits_{|\zeta|<1}
{(|\zeta|^2+1)\rho\,d\,Re\zeta\,d\,Im\,\zeta\over
|\zeta|^2(1+\rho(|\zeta|+|\zeta|^{-1}))^2|\zeta-\lambda|}=\cr
&\int\limits_{|\zeta|<1}
{\rho(|\zeta|^2+1)d\,Re\zeta\,d\,Im\,\zeta\over
(|\zeta|+\rho(|\zeta|^2+1))^2|\zeta-\lambda|},\cr}\eqno(11.8b)$$
$$J_{2,2}(\lambda,\rho)=\int\limits_{|\zeta|>1}
{(|\zeta|^2+1)\rho\,d\,Re\zeta\,d\,Im\,\zeta\over
|\zeta|^2(1+\rho(|\zeta|+|\zeta|^{-1}))^2|\zeta-\lambda|},\eqno(11.8c)$$
where $\lambda\in\C$, $\rho>0$. In addition,
$$\eqalign{
&J_{2,1}(\lambda,\rho)\le \biggl(
\int\limits_{\scriptstyle |\zeta|<1 \atop |\zeta|\le |\zeta-\lambda|}+
\int\limits_{\scriptstyle |\zeta|<1 \atop |\zeta|\ge |\zeta-\lambda|}
\biggr){2\rho\,d\,Re\zeta\,d\,Im\,\zeta\over (|\zeta|+\rho)^2|\zeta-\lambda|}
\le\cr
&\int\limits_{|\zeta|<1}
{2\rho\,d\,Re\zeta\,d\,Im\,\zeta\over (|\zeta|+\rho)^2|\zeta|}+
\int\limits_{|\zeta|<1}
{2\rho\,d\,Re\zeta\,d\,Im\,\zeta\over (|\zeta-\lambda|+\rho)^2|\zeta-\lambda|}
\le\cr
&\int\limits_{\C}
{4\rho\,d\,Re\zeta\,d\,Im\,\zeta\over (|\zeta|+\rho)^2|\zeta|}=
\int_0^{\infty}{8\pi\rho\,dr\over (r+\rho)^2}=8\pi,\cr}\eqno(11.9a)$$
$$\eqalign{
&J_{2,2}(\lambda,\rho)\le \biggl(
\int\limits_{\scriptstyle |\zeta|>1 \atop |\zeta|\le |\zeta-\lambda|}+
\int\limits_{\scriptstyle |\zeta|>1 \atop |\zeta|\ge |\zeta-\lambda|}
\biggr){2\rho\,d\,Re\zeta\,d\,Im\,\zeta\over
(1+\rho\,|\zeta|)^2|\zeta-\lambda|}
\le\cr
&\int\limits_{|\zeta|>1}
{2\rho\,d\,Re\zeta\,d\,Im\,\zeta\over (1+\rho\,|\zeta|)^2|\zeta|}+
\int\limits_{|\zeta|>1}
{2\rho\,d\,Re\zeta\,d\,Im\,\zeta\over
(1+\rho\,|\zeta-\lambda|)^2|\zeta-\lambda|}
\le\cr
&\int\limits_{\C}
{4\rho\,d\,Re\zeta\,d\,Im\,\zeta\over (1+\rho\,|\zeta|)^2|\zeta|}=
\int_0^{\infty}{8\pi\rho\,dr\over (1+\rho\,r)^2}=8\pi,\cr}\eqno(11.9b)$$
where $\lambda\in\C$, $\rho>0$. Estimate (11.5) follows from (11.8), (11.9).

\vskip 2 mm
{\it Proof of (11.6)}.
We have that
$$\eqalign{
&J_3(\lambda,\rho)\le \biggl(
\int\limits_{|\zeta|\le |\zeta-\lambda|}+
\int\limits_{|\zeta|\ge |\zeta-\lambda|}
\biggr){\rho\,d\,Re\zeta\,d\,Im\,\zeta\over
(|\zeta|+\rho(|\zeta|^2+1)|\zeta-\lambda|}
\le\cr
&\int\limits_{\C}
{2\rho\,d\,Re\zeta\,d\,Im\,\zeta\over (|\zeta|+\rho(|\zeta|^2+1))|\zeta|}=
\int_0^{\infty}{4\pi\rho\,dr\over {r+\rho(r^2+1)}}\le\cr
&\int_0^1{4\pi\rho\,dr\over {r+\rho}}+\int_1^{\infty}{4\pi\rho\,dr\over
r(1+\rho\,r)}=\cr
&\int_0^1{4\pi\rho\,dr\over {r+\rho}}=8\pi\rho\ln\,\bigl({{1+\rho}\over
\rho}\bigr),\cr}\eqno(11.10)$$
where $\lambda\in\C$, $\rho>0$. Estimate (11.6) follows from (11.10).

Lemma 11.1 is proved.

Using formulas (6.13c), (6.8), Lemmas 4.1, 5.2, 11.1 and smoothing properties
of the convolution with $1/\zeta$ on the complex plane $\C$ we obtain
properties (6.14) for $I(U_1,U_2)$ and estimate (6.15a).
Properties and estimates (6.14), (6.15b), (6.15c) for $N(U)$ and $M(U)$
follow from property (6.14) and estimate (6.15a) for $I(U_1,U_2)$. Estimate
(6.16a) follows from the formula
$$N(U_1)-N(U_2)=I(U_1-U_2,U_1)+I(U_2,U_1-U_2)   \eqno(11.11)$$
and from estimate (6.15a). Estimate (6.16b) follows from (6.13a), (6.14a)
and (6.16a).

Lemma 6.1 is proved.

\vskip 2 mm
{\bf 12. Proof of Lemmas 6.2 and 6.3}

{\it Proof of Lemma 6.2.}
Suppose that
$$U,V\in L_{\mu}^{\infty}((\C\b 0)\times (\R^3\b {\cal L}_{\nu})),\ \
|||U|||_{\mu}\le r,\ \  |||V|||_{\mu}\le r.\eqno(12.1)$$
Then using Lemma 6.1 and the assumptions of Lemma 6.2 we obtain that
$$M_{U_0}(U)\in L_{\mu}^{\infty}((\C\b 0)\times (\R^3\b {\cal L}_{\nu})),$$
$$\eqalignno{
&|||M_{U_0}(U)|||_{\mu}\le |||U_0|||_{\mu}+|||M(U)|||_{\mu}\le
r/2 + 2c_6(\mu)r^2<r,&(12.2)\cr
&|||M_{U_0}(U)-M_{U_0}(V)|||_{\mu}\le\alpha\,|||U-V|||_{\mu},\ \
\alpha=4c_6(\mu)r<1,&(12.3)\cr}$$
where
$$M_{U_0}(U)=U_0+M(U).\eqno(12.4)$$
Due to (12.1)-(12.4), $M_{U_0}$ is a contraction map of the ball
$U\in L_{\mu}^{\infty}((\C\b 0)\times (\R^3\b {\cal L}_{\nu}))$,
$|||U|||_{\mu}\le r$. Using now the lemma about contraction maps we obtain
that (6.17) is uniquely solvable for $U$ of the aforementioned ball
by the method of successive approximations. In addition, using the formulas
$$|||U-(M_{U_0})^n(0)|||_{\mu}\le \sum_{j=n}^{\infty}
|||(M_{U_0})^{j+1}(0)-(M_{U_0})^j(0)|||_{\mu},\eqno(12.5)$$
$$\eqalign{
&|||(M_{U_0})^{j+1}(0)-(M_{U_0})^j(0)|||_{\mu}\buildrel (12.3) \over \le\cr
&4c_6(\mu)r|||(M_{U_0})^j(0)-(M_{U_0})^{j-1}(0)|||_{\mu},\ \
j=1,2,3,\ldots,\cr}\eqno(12.6a)$$
$$\eqalign{
&|||(M_{U_0})^{j+1}(0)-(M_{U_0})^j(0)|||_{\mu}\buildrel (12.6a) \over \le\cr
&(4c_6(\mu)r)^j|||M_{U_0}(0)-(M_{U_0})^0(0)|||_{\mu}\buildrel (12.4) \over =
\cr
&(4c_6(\mu)r)^j|||U_0|||_{\mu}\le (4c_6(\mu)r)^j r/2,\ \
j=1,2,3,\ldots,\cr}\eqno(12.6b)$$
where $U$ is the solution of (6.17) in the  aforementioned ball and
$(M_{U_0})^0(0)=0$, we obtain (6.18).

Lemma 6.2 is proved.

\vskip 2 mm
{\it Proof of Lemma 6.3}.
We have that
$$\eqalignno{
&U-\tilde U=U_0-\tilde U_0+M(U)-M(\tilde U),&(12.7)\cr
&M(U)(\lambda,p)-M(\tilde U)(\lambda,p)\buildrel (6.13a),(11.11) \over =
L_{U,\tilde U}(U-\tilde U),&(12.8)\cr}$$
where
$$L_{U,\tilde U}W=I(W,U)(\lambda,p)+I(\tilde U,W)(\lambda,p)+
I(W,U)(\lambda_0(p),p)+I(\tilde U,W)(\lambda_0(p),p),\eqno(12.9)$$
where $I(U_1,U_2)$ is defined by (6.13c), $W$ is a test function on
$(\C\b 0)\times (\R^3\b {\cal L}_{\nu})$. In view of (12.8), (12.9) we can
consider (12.7) as a linear integral equation for "unknown" $U-\tilde U$
with given $U_0-\tilde U_0$, $U$, $\tilde U$. Using (12.9), (6.14), (6.15a),
and the properties  $|||U|||_{\mu}\le r$, $|||\tilde U|||_{\mu}\le r$,
we obtain that
$$\eqalign{
&L_{U,\tilde U}W\in L_{\mu}^{\infty}((\C\b 0)\times (\R^3\b {\cal L}_{\nu})),
\cr
&|||L_{U,\tilde U}W|||_{\mu}\le 4c_6(\mu)r|||W|||_{\mu}\ \ {\rm for}\ \
W\in L_{\mu}^{\infty}((\C\b 0)\times (\R^3\b {\cal L}_{\nu})).\cr}\eqno(12.10)
$$
Using (12.8)-(12.10) and solving (12.7) with respect to $U-\tilde U$ by the
method of successive approximations, we obtain (6.19).

Lemma 6.3 is proved.

\vskip 4 mm
{\bf References}
\vskip 2 mm
\item{[   A]} G.Alessandrini, {\it Stable determination of conductivity
by boundary measurements}, Appl. Anal. {\bf 27} (1988), 153-172.
\item{[ BC1]} R.Beals and R.R.Coifman, {\it Multidimensional inverse
scattering and nonlinear partial differential equations}, Proc. Symp. Pure
Math. {\bf 43} (1985), 45-70.
\item{[ BC2]} R.Beals and R.R.Coifman, {\it The spectral problem for the
Davey-Stewartson and Ishimori hierarchies}, Nonlinear evolution
equations: integrability and spectral methods, Proc. Workshop, Como/Italy
1988, Proc. Nonlinear Sci., (1990), 15-23.
\item{[BLMP]} M.Boiti, J.Leon, M.Manna and F.Pempinelli, {\it On a spectral
transform of a KDV- like equation related to the Schr\"odinger operator in
the plane}, Inverse Problems {\bf 3} (1987), 25-36.
\item{[  BU]} R.M.Brown and G.Uhlmann, {\it Uniqueness in the inverse
conductivity problem for nonsmooth conductivities in two dimensions},
Comm. Partial Diff. Eq. {\bf 22} (1997), 1009-1027.
\item{[   C]} A.-P.Calder\'on, {\it On an inverse boundary value problem},
Seminar on Numerical Analysis and its Applications to Continuum Physics
(Rio de Janeiro, 1980), pp.65-73, Soc. Brasil. Mat. Rio de Janeiro, 1980.
\item{[  ER]} G.Eskin and J.Ralston, {\it The inverse back-scattering problem
in three dimensions}, Commun. Math. Phys. {\bf 124} (1989), 169-215.
\item{[  F1]} L.D.Faddeev, {\it Growing solutions of the Schr\"odinger
equation}, Dokl. Akad. Nauk SSSR {\bf 165} (1965), 514-517 (in Russian);
English Transl.: Sov. Phys. Dokl. {\bf 10} (1966), 1033-1035.
\item{[  F2]} L.D.Faddeev, {\it Inverse problem of quantum scattering theory
II}, Itogi Nauki i Tekhniki, Sovr. Prob. Math. {\bf 3} (1974), 93-180
(in Russian); English Transl.: J.Sov. Math. {\bf 5} (1976), 334-396.
\item{[   G]} I.M.Gelfand, {\it Some problems of functional analysis and
algebra}, Proceedings of the International Congress of Mathematicians,
Amsterdam, 1954, pp.253-276.
\item{[  GN]} P.G.Grinevich, S.P.Novikov, {\it Two-dimensional "inverse
scattering problem" for negative energies and generalized-analytic functions.
I. Energies below the ground state}, Funkt. Anal. i Pril. {\bf 22}(1)
(1988), 23-33 (In Russian); English Transl.: Funkt. Anal. and Appl. {\bf 22}
(1988), 19-27.
\item{[  HN]} G.M.Henkin and R.G.Novikov, {\it The $\bar\pa$- equation in the
multidimensional inverse scattering problem}, Uspekhi Mat. Nauk {\bf 42(3)}
(1987), 93-152 (in Russian); English Transl.: Russ. Math. Surv. {\bf 42(3)}
(1987), 109-180.
\item{[  KV]} R.Kohn and M.Vogelius, {\it Determining conductivity by
boundary measurements II, Interior results}, Comm. Pure Appl. Math. {\bf 38}
(1985), 643-667.
\item{[   M]} N.Mandache, {\it Exponential instability in an inverse problem
for the Schr\"odinger equation}, Inverse Problems {\bf 17} (2001), 1435-1444.
\item{[ Mos]} H.E.Moses, {\it Calculation of a scattering potential from
reflection coefficients}, Phys. Rev. (2) {\bf 102} (1956), 559-567.
\item{[ Na1]} A.I.Nachman, {\it Reconstructions from boundary measurements},
Ann. Math. {\bf 128} (1988), 531-576.
\item{[ Na2]} A.I.Nachman, {\it Global uniqueness for a two-dimensional
inverse boundary value problem}, Ann, Math. {\bf 142} (1995), 71-96.
\item{[ No1]} R.G.Novikov, {\it Multidimensional inverse spectral problem
for the equation $-\Delta\psi+(v(x)-Eu(x))\psi=0$}, Funkt. Anal. i Pril.
{\bf 22(4)} (1988), 11-22 (in Russian); English Transl.: Funct. Anal. and
Appl. {\bf 22} (1988), 263-272.
\item{[ No2]} R.G.Novikov, {\it The inverse scattering problem at fixed
energy level for the two-dimen-

\item{      } sional Schr\"odinger operator}, J.Funct. Anal.
{\bf 103} (1992), 409-463.
\item{[ No3]} R.G.Novikov, {\it Scattering for the Schr\"odinger equation
in multidimensional non-linear $\bar\pa$- equation, characterization of
scattering data and related results}, Scattering

\item{      } (E.R.Pike and
P.Sabatier, eds) chapter 6.2.4, Academic, New York, 2002
\item{[ No4]} R.G.Novikov, {\it Formulae and equations for finding scattering
data from the Dirichlet-to-Neumann map with nonzero background potential},
Inverse Problems {\bf 21} (2005), 257-270.
\item{[ No5]} R.G.Novikov, {\it The $\bar\pa$- approach to approximate inverse
scattering at fixed energy in three dimensions}, International Mathematics
Research Papers, {\bf 2005:6}, (2005), 287-349.
\item{[   P]} R.T.Prosser, {\it Formal solutions of inverse scattering problem.
III}, J.Math. Phys. {\bf 21} (1980), 2648-2653.
\item{[  SU]} J.Sylvester and G.Uhlmann, {\it A global uniqueness theorem
for an inverse boundary value problem}, Ann. Math. {\bf 125} (1987),
153-169
\item{[   T]} T.Y.Tsai, {\it The Schr\"odinger operator in the plane},
Inverse Problems {\bf 9} (1993), 763-787.

\end